\documentclass{article}
\usepackage[top=2cm,bottom=2cm, left=2cm, right=3cm]{geometry}
\usepackage{amsmath, amssymb, amsthm}
\usepackage[utf8]{inputenc}

\usepackage{parskip}
\usepackage{tikz}
\usepackage{tikz}
\usepackage{tikz-cd}
\usepackage{url}
\usepackage{enumerate}
\usepackage{enumitem}
\usepackage{hyperref}
\usepackage{multicol}
\usepackage{todonotes}
\usepackage{float}

\newtheorem{Theorem}{Theorem}[section]
\newtheorem{Lemma}[Theorem]{Lemma}

\newtheorem*{Theorem*}{Theorem}

\theoremstyle{definition}
\newtheorem{Definition}[Theorem]{Definition}
\newtheorem{Remark}[Theorem]{Remark}
\newtheorem{Example}[Theorem]{Example}
\usepackage[utf8]{inputenc}
\usepackage[T1]{fontenc}
\usepackage[
backend=biber, 
style=alphabetic, 
maxnames=99,
]{biblatex}
\DeclareSourcemap{
	\maps[datatype=bibtex]{
		\map[overwrite=true]{
			\step[fieldset=urldate, null]
		}
	}
}

\title{Darboux type theorems in multisymplectic geometry}
\author{
Leonid Ryvkin\thanks{
Université Claude Bernard Lyon 1, France and University of Göttingen, Germany. \texttt{ryvkin@math.univ-lyon1.fr}
},
}

\begin{document}

\maketitle

\begin{abstract}
We give a survey of Darboux type theorems in multisymplectic geometry. These theorems establish when a closed differential form of a certain type admits a constant-coefficient expression in some local coordinate system. Beyond the classical cases of symplectic and volume forms, 0-deformability (i.e. constancy of linear type) is typically not automatic and has to be imposed, leading to distinct theorems 'per linear type'. 
\end{abstract}

\tableofcontents

\section{Introduction}
Multisymplectic manifolds were originally introduced in \cite{kijowskiFinitedimensionalCanonicalFormalism1973} to give a covariant and finite-dimensional description of classical field theory, generalizing the way how symplectic geometry describes classical-mechanics. The symplectic two-form therein is replaced by a (potentially higher degree) form $\omega$, which is closed and non-degenerate (in the sense that for tangent vectors $v$ the equation $\iota_v\omega=0$ implies $v=0$). We refer to \cite{ryvkinInvitationMultisymplecticGeometry2019} for a recent introduction to the subject. 

Given a class of multisymplectic forms it is very natural to ask for Darboux theorems. Darboux theorems assert that a given form always admits a standard (constant coefficient) expression if an appropriate coordinate system is chosen. Whether a Darboux theorem exists depends heavily on what type of form one is looking at: It is true for symplectic forms and volume forms (cf. Theorem \ref{thm:darbouxclassical}), but false in general for multisymplectic 3-forms in 5 dimensions, where an additional involutivity condition needs to be required in order to have a nice coordinate expression (cf. Theorem \ref{thm:dansymp}). In higher dimensions the situation gets even more chaotic: Already in dimension six, two constant-coefficient non-degenerate 3-forms can be inequivalent, as for instance the following forms (cf. e.g. Theorem \ref{thm:binary}):
\begin{align}\label{eq:exdifftypes}
&dx_1\wedge dx_2\wedge dx_3 + dx_4\wedge dx_5\wedge dx_6\\
&dx_1\wedge dx_4\wedge dx_5 + dx_2\wedge dx_4\wedge dx_6
+ dx_3\wedge dx_5\wedge dx_6.\nonumber
\end{align}

The goal of the present survey is to give an overview of the different types of multisymplectic forms for which criteria for 'Darboux type theorems' exist in the literature. 
It is somewhat transversal to the following recent articles:

\begin{itemize}
\item \cite{graciaDarbouxTheoremsGeometric2024}: This article contains interesting results on geometries defined by multiple differential forms (in addition to certain multisymplectic results which we will also review here).
\item \cite{leClassificationKformsRn2020}: This article covers linear normal form results for exterior forms over various fields, including certain cases where there is an infinite number of possible normal forms (while we only discuss the $\mathbb R$-case and only cases where the number of possible normal forms is finite).
\end{itemize}

Most of the material in the present work is not new (except maybe for Theorem \ref{thm:martininvolutive}), but we try to summarize various developments scattered in the literature. In case I missed something, please let me know! This article contains a modification of parts of the Phd thesis \cite{ryvkinNormalFormsConserved2018} and the unpublished notes \cite{ryvkinLinearOrbitsAlternating2017,ryvkinMultisymplecticManifoldNot2016}.

\subsection*{Acknowledgements}

The author would like to thank Udhav Fowdar, Laura Leski, Julieth Saavedra, Bernd Stratmann, Maxime Wagner and Tilmann Wurzbacher for various discussions around the topic of this survey. Special thanks go to Luca Vitagliano and Alfonso Tortorella for their encouragement to write up these notes. This research was supported by the DFG grant Higher Lie Theory - Project number 539126009. The author thanks the anonymous referee for helpful comments and corrections. For the purpose of Open Access, a CC-BY-NC-SA public copyright license has been applied by the authors to the present document and will be applied to all subsequent versions up to the Author Accepted Manuscript arising from this submission.

\subsection{Preliminaries}
We will deal with smooth objects (manifolds, differential forms, vector fields etc.) over $\mathbb R$ throughout, unless otherwise mentioned. 

When we look at a differential form $\omega\in \Omega^k(M)$ on an $n$-dimensional manifold $M$ the first type of information we can extract from it is point-wise, i.e. how does $\omega_p$ look as an exterior form in $\Lambda^kT^*_pM$. Of course, in general $T_pM$ does not have a natural preferred basis, so we should think of the linear type of $\omega_p$ as being 'up to basis change'. This motivates the following definition:

\begin{Definition} Let $n\in\mathbb N$ and $k\in\{1,...,n\}$. A linear type of $k$-forms in dimension $n$ (or shortly $(k,n)-type$) is an element of the quotient $\frac{\Lambda^k(\mathbb R^n)^*}{GL(\mathbb R^n)}$, with  respect to the natural action of $GL(\mathbb R^n)$ on ${\Lambda^k (\mathbb R^n)^*}$.\footnote{i.e. $(g\cdot \alpha)(v_1,...,v_k)=\alpha(gv_1,...,gv_k)$.} 
For any vector space $V$ (e.g. $V=T_pM$), we say that $\omega\in \Lambda^kV^*$ is of type $T\in \frac{\Lambda^k(\mathbb R^n)^*}{GL(\mathbb R^n)}$, if $\phi^*\omega$ is of type $T$ for any (or equivalently every) isomorphism $\phi:\mathbb R^n\to V$. We will usually write $T$ simply as a form, leaving the fact that we mean its equivalence class implicit.

\end{Definition}

\begin{Example} The forms from Equation \eqref{eq:exdifftypes} correspond to different $(3,6)$-types.
\end{Example}

We will now introduce the central notion for this article - flatness:

\begin{Definition}\label{def:lintyp} Let $M$ be an $n$-dimensional manifold. A differential form $\omega\in\Omega^k(M)$ is called flat at $p$ if there exist local coordinates $\phi:U\subset \mathbb R^n\to M$ near $p$, such that $\phi^*\omega$ has constant coefficients, i.e. 
	$$
	\phi^*\omega=\sum_{1\leq i_1<i_2...<i_k\leq n} \lambda_{i_1,...,i_k} dx_{i_1}\wedge ...\wedge dx_{i_k}
	$$ 	
for some constants $\lambda_{i_1,...,i_k}\in \mathbb R$. We call $\omega$ flat if it is flat near any $p\in M$.
\end{Definition}

Flatness has multiple immediate consequences:

\begin{Lemma} Let $\omega$ be flat then:
	\begin{enumerate}
		\item It has constant linear type\footnote{I.e.,	 at all points it corresponds to the same linear type $T\in \frac{\Lambda^k(\mathbb R^n)^*}{GL(\mathbb R^n)}$ as described in Definition \ref{def:lintyp}.}.
		\item It is closed.
	\end{enumerate}
\end{Lemma}

\begin{Example} ~
\begin{itemize}
	\item Let $M=\mathbb R$ The form $dx$ is flat. The form $xdx$ is not flat, because it has different linear types for $x=0$ and $x\neq 0$.
	\item Let $M=\mathbb R^2$ The form $e^xdx\wedge dy$ is flat (by changing coordinates to $\tilde x=e^x, \tilde y=y$). The form $ydx$ is not flat, since it is not closed.
\end{itemize}
\end{Example}

The goal of this article is studying the flatness of differential forms. As indicated above, we can - without loss of generality - assume the forms we consider to be closed, because otherwise they have no chance of being flat. Interestingly, the other ingredient to multisymplecticity - non-degeneracy - also automatically appears in the context of studying flatness. Let us recall the definition.

\begin{Definition} ~ 
	\begin{itemize}
	\item Let $V$ be a vector space. A form $\omega\in \Lambda^kV^*$ is called non-degenerate, if the map 
	$
	V\to \Lambda^{k-1}V^*, v\mapsto \iota_v\omega
	$
is injective.
\item Let $M$ be a smooth manifold and $\omega\in \Omega^k(M)$. The form $\omega$ is called multisymplectic and the pair $(M,\omega)$ a multisymplectic manifold if $\omega$ is closed and non-degenerate.
\end{itemize}
\end{Definition}

The following Lemma asserts that any flat form locally looks like the pullback of a non-degenerate flat form. Locally, pulling back along a submersion can be understood as introducing some additional coordinates that don't appear in the form. Hence, in the sequel we can restrict our attention to non-degenerate forms, when thinking about flatness (cf. also Remark \ref{rem:flatdis} below).

\begin{Lemma} Let $\omega\in\Omega^k(M)$ and $p\in M$. The following are equivalent:
	\begin{enumerate}
		\item $\omega$ is flat near $p$.
		\item There is an open neighborhood $U$ of $p$, a surjective submersion $\phi:U\to \bar U\subset \mathbb R^{\bar n}$ and a flat non-degenerate $\bar\omega \in \Omega^k(\bar U)$ such that $\omega=\phi^*\bar\omega$.
	\end{enumerate} 
\end{Lemma}

\begin{proof}
By the rank theorem, locally, any surjective submersion is equivalent to a projection $$(x_1,...,x_{n-r},y_1,...,y_r)\to (x_1,...,x_k).$$ With respect to these coordinates the pullback along $\phi$ just means adding the additional coordinates $y_i$ to the domain without them appearing in the form. A coordinate transformation $\psi$ turning $\bar\omega$ into a constant coefficient form can be reinterpreted as the transformation $\Psi$ with $\Psi(x,y)=(\psi(x),y)$ which will turn $\omega$ into a constant coefficient form. This already implies the direction $2.\Rightarrow 1.$.
For the converse direction, assume that $\omega$ is flat near $p$ and pick local coordinates $(\tilde x_1,...,\tilde x_n)$ on $U\subset M$ such that it has constant coefficients. Let the space  
\begin{align*}
	K=K(\omega_p):=\{v\in T_pM ~|~\iota_v\omega_p=0\}
\end{align*} 
be $r$-dimensional. Since $\omega$ has constant coefficients $K$ lies in the kernel of $\omega_q$ for all $q\in U$. Now let $A$ be an invertible matrix whose last $r$ columns span $K$. We now carry out the linear coordinate change 

$$
(\tilde x_1,...,\tilde x_n)=A(x_1,...,x_{n-r},y_1,...,y_r).
$$
With respect to these new coordinates, $K$ is generated by $\partial_{y_1},...,\partial_{y_r}$, i.e. $\omega$ is a form only in the $x$ coordinates. Since the coordinate transformation was linear, $\omega$ has constant coefficients also with respect to the new coordinates, showing the direction $1.\Rightarrow 2.$.  

\end{proof}

\begin{Remark}\label{rem:flatdis} Even without requiring flatness, we can consider the space $K(\omega):=\{v\in TM ~|~\iota_v\omega=0\}$ for a differential form. Its $M$-fibers are vector spaces (possibly with dimensions jumping upward at certain points). It satisfies: 
\begin{enumerate}
		\item If $K(\omega)$ has a constant dimension, then it forms a regular distribution.\\
		 This is true because $K(\omega)$ is the kernel of the (constant rank) vector bundle morphism $TM\to \Lambda^{k-1}T^*M, v\mapsto \iota_v\omega$.
		\item If $\omega$ is closed, $\Gamma(K(\omega))$ is involutive.\\
		This can be shown as follows: Let $X,Y\in \Gamma(K(\omega))$. We want to show that their Lie bracket also is in $K(\omega)$. We verify this using Cartan calculus:
		\begin{align*}
			\iota_{[X,Y]}\omega=L_X\iota_Y\omega -\iota_YL_X\omega=-\iota_YL_X\omega=-\iota_Yd\iota_X\omega-\iota_Y\iota_Xd\omega=0.
		\end{align*}
	\end{enumerate}
When $\omega$ is closed and $K(\omega)$ has constant dimension this means that $K(\omega)$ forms a foliation, which allows us to (locally) see $\omega$ as the pullback of a non-degenerate form. This is what happened in the $1.\Rightarrow 2.$ direction in the above proof (under the simplifying assumption of $\omega$ having a constant coefficient representation).
\end{Remark}

\subsection{Summary}

The kind of result we want to survey are Darboux type theorems, i.e. theorems that, given a linear type (or some equivalent information), tell us when a (multisymplectic) form of this type is flat. To do so, we will first give an overview over some results on linear types which can be summarized as follows (Theorem \ref{thm:lin} in the text):

\begin{Theorem*}Let $\Sigma^k_n$ denote the number of non-degenerate linear $(k,n)$-types.
	\begin{itemize}
		\item $\Sigma_n^n=1$ for all $n$, and $\Sigma^1_n$ as well as $\Sigma^{n{-}1}_n$ are zero for $n>1$.
		\item $\Sigma^2_n$ is 0 for $n$ odd and one for $n$ even.
		\item $\Sigma_{n}^{n-2}=\lfloor \frac{n}{2}\rfloor-1$, when $(n \mod 4)\neq 2$ (for $n\geq 4$) and
		$\Sigma_{n}^{n-2}= \frac{n}{2}$, when $(n \mod 4)=2$ (for $n\geq 4$).
		\item $\Sigma_6^3=3$, $\Sigma_7^3=8$  ,  $\Sigma_8^3=21$,  $\Sigma_7^4=15$ and $\Sigma_8^5=31$.
		\item $\Sigma_n^k=\infty$ in all other cases. These infinities are uncountable.
	\end{itemize}
\end{Theorem*}

The different parts of the theorem are due to \cite{gourewitchLalgebreTrivecteur1935,martinetSingularitesFormesDifferentielles1970,westwickRealTrivectorsRank1981,djokovicClassificationTrivectorsEightdimensional1983,ryvkinLinearOrbitsAlternating2017}, we will comment on the individual contributions when we discuss the parts of the theorem in Section \ref{sec:lin}. A more extensive discussion, containing results over different fields and their relations can be found in \cite{leClassificationKformsRn2020}. I would also like to point out the recent article \cite{borovoiClassificationRealTrivectors2022}, where the (infinitely many) $(3,9)$-types were classified.\\

Then, in Section 3, we will turn to Darboux type theorems. Here is an overview of the types of forms we will look at.

\begin{itemize}
	\item $(n,n)$: $dx_1\wedge ... \wedge dx_n$ (volume forms), Theorem \ref{thm:darbouxclassical}.
	\item $(2,2m)$: $dx_1\wedge dx_2 + ...+ dx_{2m-1}\wedge dx_{2m}$ (symplectic forms), Theorem \ref{thm:darbouxclassical}.
	\item $(k, {m\choose k-1} + m)$ for $m>k>1$: $\sum_{i_1<...<i_k}dp_{i_1,...,i_k}\wedge dq_{i_1}\wedge ... \wedge dq_{i_k}$ (multicotangent type forms), Theorems \ref{thm:martin}, \ref{thm:martininvolutive}.
	\item $(m,d\cdot m)$, $m>2,d>1$: $(dx_1\wedge...\wedge dx_m)+(dx_{m+1}\wedge...\wedge dx_{2m})+...+(dx_{{(d-1)}\cdot m+1}\wedge...\wedge dx_{d\cdot m})$ \\ (product type forms), Theorem \ref{thm:prod}.
	\item $(m,2m)$: $\mathrm{Re}\left( (dx_1+idy_1)\wedge ...\wedge (dx_m+idy_m) \right)$ 
	(real part of a complex volume),  Theorem \ref{thm:binary}.
	\item $(d+2,d+2m)$ with $m>1$:  $(dx_1\wedge dx_2+.... +dx_{2m-1}\wedge dx_{2m})\wedge dy_1\wedge ...\wedge dy_d$\\ (density-valued symplectic forms), Theorem \ref{thm:dansymp}.
	\item $(n-2,n)$ with $n>4$: $ \pm \left( \sum_{i=1}^m dx_{2i-1}\wedge dx_{2i}\right)^{m-1}\wedge dy_1\wedge ... \wedge dy_r$ for $n=2m+r$. \\(forms of codegree two, Theorem \ref{thm:codeg}).
\end{itemize}
Only the first two of these types have no additional involutivity condition in general. However, when one is in a lucky rank/dimension combination, unconditional Darboux theorems can also occur for some of the others. At the end of Section \ref{sec:flat}, we note on a few further classes of multisymplectic manifolds which have constant linear type, but are not flat in general.\\

\section{Linear types of multisymplectic manifolds}
\label{sec:lin}

In this section we are going to look at the number of linear $(k,n)$-types, depending on $k$ and $n$. We are going to formulate certain phenomena in terms of a vector space $V$ rather than $\mathbb R^n$, to avoid getting confused when we use both $V$ and $V^*$. Recall that by a linear type we mean an orbit of the natural $GL(V)$-action on $\Lambda^kV^*$.

\begin{Remark}\label{rem:dual}
Some of the literature, e.g.
\cite{westwickRealTrivectorsRank1981, djokovicClassificationTrivectorsEightdimensional1983}, works with multi-vectors (i.e. elements of $\Lambda^kV$) rather than forms, however these two perspectives are equivalent. Any isomorphism between $V$ and $V^*$ induces a bijection on the level of $GL(V)$-orbits.
\end{Remark}

We are going to focus on non-degenerate types, since they retain all the information we need: 

\begin{Lemma}\label{degenerateLemma}
	The degenerate linear $(k,n)$-types are in bijection with all linear $(k,n-1)$-types.
\end{Lemma}
\begin{proof}
	Let $W$ be a $(n{-}1)$-dimensional subspace of $V$ and $a\in V\backslash W$. Then $V=W\oplus \mathbb R\cdot a$. Let $\text{pr}_W:V\to W$ be the projection. We show that the $GL(W)$-orbits of $\Lambda^k(W^*)$ are in one-to-one correspondence to $GL(V)$-orbits of $\Lambda^k_{deg}(V^*)$, the degenerate $k$-forms on $V$, via the map 
	\[\bar \phi:\frac{\Lambda^k(W^*)}{GL(W)}\to \frac{\Lambda^k_{deg}(V^*)}{GL(V)}, ~ [\alpha]=GL(W)\cdot \alpha\mapsto GL(V)\cdot \phi(\alpha)=[\phi(\alpha)],\]
	induced by 
	\[\phi:\Lambda^k(W^*)\to \Lambda^k(W^*)\oplus(\Lambda^{k-1}(W^*)\otimes ( \mathbb R\cdot a)^*)\cong\Lambda^k(V^*),~~\alpha\mapsto (\alpha \oplus 0).\] 
	
	{\bf Surjectivity of $\bar\phi$.} Assume $\alpha\in \Lambda^k(V^*)$ is degenerate, i.e. there exists a $v$ in $V\backslash \{0\}$ such that $\iota_v\alpha=0$. Choosing a $g\in GL(V)$ such that $g^{-1}v=a$, we get $\iota_a(g\cdot \alpha)=0$, i.e. $(g\cdot \alpha)|_W\in \Lambda^k(W^*)$ satisfies $\phi((g\cdot \alpha)|_W)=g\cdot \alpha$. Consequently $\bar\phi([(g\cdot \alpha)|_W])=[\alpha]$ holds and the map $\bar \phi$ is surjective.
	
	{\bf Injectivity of $\bar\phi$.} Let $\alpha,\beta\in \Lambda^k(W^*)$ and $\phi(\alpha)$ be equivalent to $\phi(\beta)$, i.e. there exists a map $g\in GL(V)$ such that $g\cdot \phi(\alpha)=\phi(\beta)$. Especially $ga$ and $a$ are both elements of the vector subspace $\text{ann}(\phi(\alpha)):=\{v\in V| \iota_v\phi(\alpha)=0\}$. We pick an element $h\in GL(V)$ such that $h(ga)=a$ and $h\cdot \phi(\alpha) =\phi(\alpha)$.
	Then we have $(hg)\cdot\phi(\alpha)=g\cdot(h\cdot\phi(\alpha))=g\cdot \phi(\alpha)=\phi(\beta)$, where $hg(a)=a$. Consequently, $\text{pr}_W\circ (hg)|_W$ yields a well-defined automorphism of $W$ satisfying $(\text{pr}_W\circ (hg)|_W)\cdot \alpha=\beta$. Hence, $[\alpha]=[\beta]$ and $\bar\phi$ is injective.
\end{proof}

\begin{Remark} The rank of a form $\alpha$ on $\mathbb R^n$ is $n-c$, where $c$ is the dimension of $\{v|\iota_v\alpha=0\}$. This rank coincides with the dimension in which $\alpha$ would define a non-degenerate linear type. Non-degeneracy means exactly that the rank is full. 
\end{Remark}

In addition to linear types and non-degenerate linear types, we are going to note which of the types we discuss are stable:

\begin{Definition} The linear type of a form $\omega\in \Lambda^kV^*$ is called stable if $GL(V)\cdot \omega$ is open in $\Lambda^kV^*$. 
\end{Definition}

\begin{Remark}\label{rem:rank} Whenever non-degenerate types exist all stable types are non-degenerate. This is because, when they exist, non-degenerate forms form an open and dense subset of all forms. To see this, pick a scalar product and write for any $\beta \in\Lambda^kV^*$ $L_\beta$ for the map $V\to \Lambda^{k-1}V^*$, $v\mapsto L_\beta(v)=\iota_v\beta$. The expression $P(\beta)=det(L_\beta^*\circ L_\beta)$, where $(-)^*$ denotes the adjoint, is a polynomial in $\beta$. This polynomial is non-zero exactly when $\beta$ is non-degenerate. Since $P$ is not constantly zero, it is non-zero on an open and dense subset, i.e. no degenerate form has a degenerate open neighborhood.
\end{Remark}

Let us exemplify the above notions by looking at the case of $(n,n)$-types:

\begin{Lemma}\label{vollemma} There are two linear $(n,n)$-types. With respect to a basis $e^1,...,e^n$ of $V^*$ they are represented by:
\begin{enumerate}
\item $e^1\wedge e^2\wedge ... \wedge e^n$ (stable, non-degenerate)
\item $0$ (not stable, degenerate)
\end{enumerate}
\end{Lemma}

\begin{proof} This follows from the fact that $GL(V)$ acts on the one-dimensional $\Lambda^nV^*$ via multiplication with the determinant.
\end{proof}

In the sequel, we are going to look at the classification of the remaining cases in three pieces:
\begin{itemize}
	\item The generic case (where $n$ is sufficiently big and $k$ is not in $\{1,2,n-2,n-1\}$). 
	\item The remaining cases with small $k$, ($k\in \{1,2\}$, and $k=3$ when $n\in\{6,7,8\})$.
	\item The remaining cases with large $k$
	(i.e. $k\in \{n-2,n-1\}$ and $k=n-3$ when $n\in\{7,8\})$
\end{itemize}
The first piece has been solved by \cite{martinetSingularitesFormesDifferentielles1970}, the second one is classical for $k\in\{1,2\}$ and has been resolved by \cite{gourewitchLalgebreTrivecteur1935,westwickRealTrivectorsRank1981, capdevielleClassificationFormesTrilineaires1973,djokovicClassificationTrivectorsEightdimensional1983} for the $k=3$ cases. The big $k$ cases follow from a duality technique presented in \cite{martinetSingularitesFormesDifferentielles1970} for the $k=n-2$ case, and were applied to the remaining cases by the author of the present in \cite{ryvkinLinearOrbitsAlternating2017}.

\subsection{The generic case}
For an $n$-dimensional vector space, the space $GL(V)$ is $n^2$ dimensional (as a smooth  manifold) and the vector space $\Lambda^kV^*$ has dimension $n\choose k$. For $n\geq 9$ and $k\in\{3,...,n-3\}$ and $(k,n)=(4,8)$ we have ${n\choose k}>n^2$. This has multiple consequences to linear types:

\begin{Lemma}[\cite{martinetSingularitesFormesDifferentielles1970}] Let $n\geq 9$ and $k\in\{3,...,n-3\}$ or $(k,n)=(4,8)$. Then:
	\begin{itemize}
		\item The number of non-degenerate linear $(k,n)$-types is (uncountably) infinite.
		\item There are no stable linear $(k,n)$-types.
	\end{itemize}
\end{Lemma}

\begin{proof}Since the dimension of the Lie group acting is less than the dimension of the space it is acting on, no open orbits can exist and the number of orbits is (uncountably) infinite. The only thing we need to show is that the number of non-degenerate types is infinite as well. Since being non-degenerate is an open condition in $\Lambda^kV^*$, we only need to show that one non-degenerate form exists. We will show it inductively as follows. For $n\in\{5,6,7\}$ the following are non-degenerate 3-forms (where $\{e^1,...,e^n\}$ are a basis of $V^*$):
\begin{align*}
 &\omega_5=(e^1\wedge e^{2}+e^3\wedge e^4)\wedge e^5,\\
 &\omega_6=e^1\wedge e^2\wedge e^3+e^4\wedge e^5\wedge e^6\\
 &\omega_7=(e^1\wedge e^{2}+e^3\wedge e^4+e^5\wedge e^6)\wedge e^7.
\end{align*}
Given a non-degenerate 3-form $\omega_n$ in dimension $n$ we can obtain one in dimension $n+3$ by $\omega_{n+3}=\omega_n+e^{n+1}\wedge e^{n+2}\wedge e^{n+3}$. This means that non-degenerate 3-forms $\omega_n$ exist for all $n\geq 5$. Now, given $(k,n)$ we can construct a non-degenerate $k$-form in dimension $n$ as:
\begin{align*}
\omega_{n-(k-3)}\wedge e^{n-(k-3)+1}\wedge ...\wedge e^{n}.
\end{align*}
\end{proof}

\subsection{Small $k$ cases}
By cases with small $k$ we mean cases where $k$ is not in the 'generic' interval described in the previous section, but takes smaller values. This means the cases where $k\in\{1,2\}$ or $k=3$ when $n\in\{6,7,8\}$.\footnote{The case $(k,n)=(3,5)$ will be treated in subsection \ref{sec:bigk}}

Let us start with the $k=1$ case:
\begin{Lemma}\label{lem:1n} For $n>1$ there are two (1,n)-types: $e^1$ and 0. Both are degenerate, however the former type is stable.
\end{Lemma}

For the $k=2$ case the situation is already a little more complicated and boils down to the symplectic basis theorem:

\begin{Lemma}[Symplectic basis theorem]\label{lem:2n}There are $\lfloor \frac{n}{2}\rfloor+1$ $(2,n)$-types. They are represented by
\begin{align*}
	\omega_r=\left( \sum_{i=1}^{r} e^{2i-1}\wedge e^{2i} \right)
\end{align*}
for $r \in \{0,...,\lfloor \frac{n}{2}\rfloor\}$. If $n$ is odd, no orbit is non-degenerate otherwise only $\omega_{\lfloor \frac{n}{2}\rfloor}$ is.  In either case  $\omega_{\lfloor \frac{n}{2}\rfloor}$ is the unique stable orbit.
\end{Lemma}

\begin{proof}
	Let $\alpha\in\Lambda^2V^*$ be given. Instead of finding an element $g\in GL(V)$, which transforms $\alpha$ into one of the orbits representatives, we will construct a basis $\{v_1,...,v_n\}$ such that $\alpha$ already looks like one of the above forms with respect to the dual basis $\{v^1,...,v^n\}$. The statements then will follow, because $GL(V)$ acts transitively on the space of bases of $V$.\\
	We will first assume $\alpha$ is non-degenerate and proceed inductively. If $V$ is zero-dimensional there is nothing to prove. So assume $dim(V)>0$. Let $v_n\neq 0$. As $\alpha$ is non-degenerate there exists some $\tilde v_{n-1}$ in $V$ such that $\alpha(\tilde v_{n-1},v_n)\neq 0$. Then we set $v_{n-1}=\frac{1}{\alpha(\tilde v_{n-1},v_n)}\tilde v_{n-1}$. Next we regard the kernel $\tilde V$ of the surjective linear map
	\[
	V\to\mathbb R^2,~~~v\mapsto {{\alpha(v,v_{n-1})}\choose {\alpha(v,v_n)}}
	\]
	By construction $\tilde V$ is an $n-2$-dimensional vector space on which $\alpha$ is non-degenerate, so we can repeat the above procedure until we arrive at a basis $\{v_1,...,v_n\}$  of $V$. Then $\alpha=v^1\wedge v^2+....+v^{n-1}\wedge v^n$ by construction. Especially $n$ is even. The basis for degenerate $\alpha$ can be retrieved from the non-degenerate case together with the considerations of  Lemma \ref{degenerateLemma}.\\
	Now it remains to show that the orbit of $ \sum_{i=0}^{\lfloor \frac{n}{2}\rfloor} e^{2i}\wedge e^{2i+1}$ is stable. To see this we reinterpret $\Lambda^nV^*$ as skew-symmetric matrices. Then the orbit of  $\omega_{\lfloor \frac{n}{2}\rfloor}$ is the orbit of matrices with maximal possible rank (full rank when $n$ is even, else rank $n-1$). Since the rank of a matrix is lower semicontinuous the matrices with maximal rank form an open subset. The others don't contain any open subset, because $\omega_{r}+\epsilon\omega_{\lfloor \frac{n}{2}\rfloor}$ is of rank $2\cdot \lfloor \frac{n}{2}\rfloor$ for arbitrarily small $\epsilon$.
\end{proof}

Finally, we will cite the results for $k=3$ and $n\in\{6,7,8\}$:
\begin{Lemma}[\cite{gourewitchLalgebreTrivecteur1935,westwickRealTrivectorsRank1981,djokovicClassificationTrivectorsEightdimensional1983}]~\label{lem:3insmall} \begin{itemize}
\item There are 6 linear $(3,6)$-types, 3 of them non-degenerate, two of the non-degenerate ones are stable.
\item There are 14 linear $(3,7)$-types, 8 of them non-degenerate, two of the non-degenerate ones are stable.
\item There are 35 linear $(3,8)$-types, 21 of them non-degenerate, three of the non-degenerate ones are stable.
\end{itemize}
\end{Lemma}

A list of the non-degenerate types can be found in Appendix \ref{app:threeclass}.  We could use Lemma \ref{degenerateLemma} to obtain all the types from them. However, we would need the additional information about all linear $(3,5)$-types. Fixing a basis $\{e^1,...,e^5\}$, these are:
	\begin{itemize}
		\item 0
		\item $e^1\wedge e^2\wedge e^3$
		\item $(e^1\wedge e^2+e^3\wedge e^4)\wedge e^5$
	\end{itemize}
This could be obtained using a by hand calculation, but also will follow from Lemma \ref{n-2-lemma} in the next section.

\subsection{Big $k$ cases}
\label{sec:bigk}
In this section, we will use a technique presented in \cite{martinetSingularitesFormesDifferentielles1970} to deduce the number of $(k,n)$-types (with big $k$) from the corresponding $(n-k,n)$-types.

Given a volume form $\Omega\in \Lambda^nV^*\backslash \{0\}$, we can construct a linear isomorphim $L:=\iota_\bullet\Omega:\Lambda^kV\to \Lambda^{n-k}V^*$ given by $\eta\mapsto \iota_\eta\Omega$. As it turns out, this linear isomorphism does not necessarily induce a bijection between the $\Lambda^kV$- and the $\Lambda^{n-k}V^*$-orbits. A multivector $\eta \in\Lambda^kV$ need not be in the same $GL(V)$-orbit as $(-\eta)$, when $k$ is even. (When $k$ is odd we have $(-id)_*\eta=-\eta$\footnote{where we use the lower star to indicate the action on exterior powers of $V$ and the higher star for the action on exterior powers of $V^*$}). As we will see below, this is a consequence of the fact that we can not take even roots of negative numbers, as we are in the algebraically non-closed field $\mathbb R$. The following lemma helps specify for which $(n,k)$ this phenomenon may occur.

\begin{Lemma}[\cite{martinetSingularitesFormesDifferentielles1970}] \label{lem:splitorb}	Let $V$ be $n$-dimensional and $\Omega\in \Lambda^nV^*\backslash \{0\}$ a volume form. Then the linear isomorphism $L=\iota_\bullet\Omega$ satisfies: 
	\begin{enumerate}
		\item For any $\eta$, $L$ identifies the union of the orbits of $\eta$ and $-\eta$ with the union of the orbits of $L(\eta)$ and $-L(\eta)$.
	 	\item Let $\eta$ be an element of $\Lambda^kV$. Then $L(\eta)$ and $-L(\eta)$ lie in the same $GL(V)$-orbit if and only if there exists a $g\in GL(V)$ satisfying $g_* \eta=-\det (g)\eta$.  A sufficient condition for this is for the stabilizer of $\eta$ to contain an element of negative determinant.
		\item When $n-k$ is odd, $L$ induces a surjection from $(k,n)$-types to $(n-k,n)$-types.
		\item When $k$ is odd, $L^{-1}$ induces a surjection from $(n-k,n)$-types to $(k,n)$-types.
	\end{enumerate}

\end{Lemma}

\begin{proof}
	The parts of the Lemma are a consequence of the equation
\begin{align}\label{eq:dualitycentral}
g^* L(\eta)=g^*(\iota_{\eta}\Omega)=\iota_{(g^{-1})_*\eta}g^*\Omega=det(g) \iota_{g^{-1}_*\eta}\Omega=det(g)\cdot L(g^{-1}_*\eta).
\end{align}

	More precisely:
	\begin{enumerate}
		\item We show that the image of the orbit of $\eta$ is contained in the union of the orbits of $L(\eta)$ and $-L(\eta)$. Let $\tilde \eta=g_*\eta$. We apply Equation \eqref{eq:dualitycentral}:
		$$
		L(\tilde\eta)=L(g_*\eta)=\det(g)g^{-1}L(\eta)
		$$ 
		We now set $h=\sqrt[n-k]{|\det(g)|}\cdot g^{-1}$, then: 
\begin{align}\label{eq:duality2ndeq}
	L(\tilde\eta)=sgn(\det(g))h^*L(\eta)=h^*L(sgn(\det(g))\eta)
\end{align}
 where $sgn\in\{+1,-1\}$ denotes the sign. The opposite direction works analogously.
 \item By Equation \eqref{eq:dualitycentral} the equality $g^*L(\eta)=L(-\eta)$ is equivalent to 
 $-\eta=det(g)\cdot g^{-1}_*\eta$, which in turn translates to the expression in question. An element in the stabilizer of $\eta$ with negative determinant satisfies this equation.
 \item When $n-k$ is odd, then $(-id)^*L(\eta)=-L(\eta)=L(-\eta)$, i.e. the orbits of $L(\eta)$ and $L(-\eta)$ always coincide.
 \item This is proven analogously to the previous item.

\end{enumerate}	
\end{proof}

\begin{Remark}\label{rem:dual2} Since the orbits in $\Lambda^kV$ are in one-to-one correspondence with those in $\Lambda^kV^*$ (cf. Remark \ref{rem:dual}), we can use the above Lemma to relate $(k,n)$-types with $(n-k,n)$-types. To talk about non-degeneracy of elements in $\Lambda^kV$ we  interpret them as forms on $V^*$ (i.e. elements of $\Lambda^k(V^* )^*$). In this context the following are true:
\begin{enumerate}
	\item Degenerate types $\eta$ always contain an element of negative determinant in their stabilizer.
	\item $L(\eta)$ is non-degenerate if and only if $\eta=v\wedge \tilde \eta$ for some $\tilde \eta \in\Lambda^{k-1}V$. 
\end{enumerate}
To see the first item, one simply takes a nonzero vector $e_1$ with $\iota_{e_1}\eta=0$ and extends it to a basis $e_1,...,e_n$ Then $g\in GL(V)$ defined by $g(e_1)=-e_1$ and $g(e_j)=e_j$ for $j\neq 1$ stabilizes $\eta$ and has negative determinant. The second one can be seen by choosing a basis starting with $v$ (or respectively stating with an element in the kernel of $L(\eta)$).
\end{Remark}

Let us start with the corank-one case:
\begin{Lemma} For $n>1$ there are two (n-1,n)-types: $e^2\wedge ... \wedge e^n$ and 0. Both are degenerate, however the former type is stable.
\end{Lemma}

\begin{proof}
One can show this directly, but we will try to use Lemma \ref{lem:splitorb} as a preparation and practice. Since $k$ is odd, the only thing we need to check is that none of the $(1,n)$-types from Lemma \ref{lem:1n} is split into two, i.e. that both of them have an element of negative determinant in their stabilizer. Since both are degenerate, this follows from Remark \ref{rem:dual2}.
\end{proof}

\begin{Lemma}[\cite{martinetSingularitesFormesDifferentielles1970}]\label{n-2-lemma}Let $dim(V)=n>4$ and Let $\{e_1,...,e_n\}$ be a basis of $V$. Let $\{e^1,...,e^n\}$ be its dual basis and $\Omega=e^1\wedge ... \wedge e^n$. If $n \not\equiv 2 \mod 4$, then there are $\lfloor \frac{n}{2}\rfloor+1$. $(2,n)$-types. They are represented by
	\begin{align*}
		\omega_r=\iota_{\left( \sum_{i=1}^{r} e_{2i-1}\wedge e_{2i} \right)}\Omega
	\end{align*}
	for $r \in \{0,...,\lfloor \frac{n}{2}\rfloor\}$.
	If $n \equiv 2 \mod 4$, then there is an additional type $\omega_{\lfloor\frac{n}{2}\rfloor}^-:=-\omega_{\lfloor\frac{n}{2}\rfloor}$. In both cases all types except $0$ and $(\iota_{e_1\wedge e_2}\Omega)$ are non-degenerate. In the former case $\omega_{\lfloor\frac{n}{2}\rfloor}$ is the unique stable type, in the latter $\omega_{\lfloor\frac{n}{2}\rfloor}$ and $\omega_{\lfloor\frac{n}{2}\rfloor}^-$ are stable.
\end{Lemma} 

\begin{proof}
	We first observe that using the classification in Lemma \ref{lem:2n} one can see that $\eta$ is in the same orbit as $-\eta$ for all elements of $\Lambda^2V$. Moreover, by Remark \ref{rem:dual2}, $\iota_\eta\Omega$ is equivalent to $-\iota_\eta\Omega$ for all degenerate $\eta$. Bivectors of full rank appear only when $n$ is even. Then they are given by $\eta=\sum_{i=1}^{n/2} e_{2i-1}\wedge e_{2i} $. 
	
	When $n/2$ is even, we construct a linear isomorphism $g$ satisfying the condition from Lemma \ref{lem:splitorb} item 2., as the following:
	\[
	g(e_i)=\left\{ \begin{matrix}
		e_{i-1}&i\text{ even}\\
		e_{i+1}&i\text{ odd}
	\end{matrix}\right. .
	\]
	We observe that $g_*\eta=-\eta$ and $det(g)=1$ because $n/2$ is even. Hence, when $n/2$ is even $\iota_\eta\Omega$ is in the same orbit as $-\iota_\eta\Omega$.
	
	Let us now turn to the case $n/2$ is odd. Assume there exists a $g$ satisfying the equation $g_* \eta=-\det(g)\eta$. Then we regard the n-fold exterior power of $g_*\eta$:
	\[
	det(g)\eta^{n/2}=(g_* \eta)^{n/2}=(-\det(g)\eta)^{n/2}=- det(g)^{n/2}\eta^{n/2}
	\]
	But $det(g)=- det(g)^{n/2}$ is a contradiction for $n/2$ odd. Thus no such $g$ can exist and the orbits of $\iota_\eta\Omega$ and $-\iota_\eta\Omega$ are distinct.
	
	As stability is preserved by $\iota_\bullet \Omega:\Lambda^2V\to \Lambda^{n-2}V^*$, what remains to show is the non-degeneracy statement. By Remark \ref{rem:dual2}, $\iota_\eta\Omega$ is non-degenerate if and only if $\eta=v\wedge \tilde \eta$ for some $v,\eta\in V$. If $v$ or $\tilde \eta$ is zero, then $\iota_\eta\Omega=0$, otherwise it lies in the orbit of $\iota_{e_1\wedge e_2}\Omega$.
\end{proof}

\begin{Lemma}[\cite{ryvkinLinearOrbitsAlternating2017}]\label{74} There are 20 linear $(4,7)$-types, 15 of those are non-degenerate and of four the non-degenerate ones are stable.
\end{Lemma}
\begin{proof}
	Since $k=3$ is odd, we only need to check which of the 14 types decreed by Lemma \ref{lem:3insmall} split into two. By Lemma \ref{lem:splitorb} this happens to a type $\eta$ whether it does not admit an element $g$ satisfying 
	\begin{align}\label{eq:eqtosplit}
		g_*\eta=-det(g)\eta.
	\end{align} 
	This equation is satisfied by any element with negative determinant in the stabilizer of $g$. Since $k$ is odd, the converse is also true:  Imagine $g$ satisfies Equation \eqref{eq:eqtosplit}. Then $$h=g\cdot \frac{-1}{\sqrt[3]{\det(g)}}$$ stabilizes $\eta$ and has negative determinant.

	 We know that the stabilizers of  degenerate 3-forms contain elements of negative determinant (cf. Remark \ref{rem:dual2}) and by Appendix \ref{a37} we know that only two of the stabilizers of non-degenerate 3-forms contain elements of negative determinant (none of those stable). Hence, there are $2\cdot 6+8=20$ orbits of 4-forms in seven-dimensional space, 4 of them stable.  What remains to calculate is the number of non-degenerate 4-forms. By Lemmas \ref{degenerateLemma} and \ref{n-2-lemma} there are five degenerate orbits, i.e. $20-5=15$ non-degenerate ones. 
\end{proof}

\begin{Lemma}[\cite{ryvkinLinearOrbitsAlternating2017}]\label{85}There are 35 $(5,8)$-types, 31 of those are non-degenerate and three of the non-degenerate ones are stable.
\end{Lemma}
\begin{proof}
	We fix a volume form $\Omega$ on $V$. In the case at hand, $(k,n)=(3,8)$, both $k$ and $n{-}k$ are odd. Hence $L$ induces a bijection between the linear $(3,8)$ and the $(5,8)$-types. This bijection preserves stability. Thus the total number of types and the number of stable types correspond to the numbers known from the case $(k,n)=(3,8)$. To get the number of non-degenerate types one observes that in the case $(k,n)=(5,7)$ there are 4 linear types and applies Lemma \ref{degenerateLemma}.
\end{proof}

\begin{Remark} The lists of all linear types in Appendix \ref{app:threeclass} can be used to obtain lists of all linear types of $(4,7)$- and $(5,8)$-types. In the case of $(5,8)$-types this can be simply done by applying the Hodge $*$-operation to all $(3,8)$-types, i.e. replacing $e^{ijk}$ with $*e^{ijk}=\iota_{e_i\wedge e_j\wedge e_k}\Omega$ where $\Omega = e^1 \wedge ...\wedge e^8$. In the $(4,7)$-types one proceeds analogously, but obtains as distinct types $(*\eta)$ and $(-{*}\eta)$ for $\eta$ in the types marked with $(+)$ in the list provided in Appendix \ref{a37}.
\end{Remark}

\subsection{Summary}

The above results can be summarized by the following: 
\begin{Theorem}\label{thm:lin} The following are the numbers of $(k,n)$-types depending on form degree $k$ and dimension $n$. The infinities involved are uncountable.
	\begin{table}[H]
		\centering
		\begin{tabular}{|l|l|l|l|}
			\hline
			& linear types      & non-deg. linear types& stable linear types \\ \hline
			$k=1$                        & 2                           & 0                     &   1             \\ \hline
			$k=2$, $n$ even                        & $\frac{n}{2}+1$ & 1  &   1  \\ \hline
			$k=2$, $n$ odd                        & $\frac{n+1}{2}$ & 0 &   1  \\ \hline
			$k=3$, $n=6$                   & 6                           & 3        &  2                           \\ \hline
			$k=3$, $n=7$                   & 14                          & 8         &    2                        \\ \hline
			$k=3$, $n=8$                   & 35                          & 21         &     3                      \\ \hline
			$k=4$, $n=7$                   &   {20}      &{15}       & 4            \\ \hline
			$k=4, n=8$                    & $\infty$                    & $\infty$                    &  0        \\ \hline
			$k=5$, $n=8$                  &   {35}      &{31}            &    3    \\ \hline
			$3\leq k\leq n-3$, $n\geq 9$ & $\infty$                    & $\infty$                       &     0  \\ \hline
			$k=n-2$, $n \geq 6$, $n=2 \mod 4$   &   $\frac{n}{2}+2$             &      $\frac{n}{2}$    & 2                       \\\hline
			$k=n-2$, $n \geq 5$, $n\neq 2 \mod 4$   &   $\lfloor\frac{n}{2}\rfloor+1$             &       $\lfloor\frac{n}{2}\rfloor - 1$      &      1                \\ \hline
			$k=n-1$                      & 2                           & 0                      &  1             \\ \hline
			$k=n$                        & 2                           & 1                    &        1         \\ \hline
		\end{tabular}
	\end{table}
\end{Theorem}

For dimensions up to 10 the non-degenerate type numbers look as follows, where the rows range from 0-forms (the ``$-$'' in the table) to $n$-forms:\\

\begin{align*}
	\begin{array}{c|ccccccccccccccccccccc}
		n=0&&&&&&&&&&&-\\[0.5em]
		n=1&&&&&&&&&&-&&1\\[0.5em]
		n=2&&&&&&&&&-&&0&&1\\[0.5em]
		n=3&&&&&&&&-&&0&&0&&1\\[0.5em]
		n=4&&&&&&&-&&0&&1&&0&&1\\[0.5em]
		n=5&&&&&&-&&0&&0&&1&&0&&1\\[0.5em]
		n=6&&&&&-&&0&&1&&3&& 3&&0&&1\\[0.5em]
		n=7&&&&-&&0&&0&& 8&& 15&&2&&0&&1\\[0.5em]
		n=8&&&-&&0&&1&& 21&&\infty&&31&&3&&0&&1\\[0.5em]
		n=9&&-&&0&&0&&\infty&&\infty&&\infty&&\infty&&3&&0&&1\\[0.5em]
		n=10&-&&0&&1&&\infty&&\infty&&\infty&&\infty&&\infty&& 5&&0&&1\\[0.5em]
	\end{array}
\end{align*}

\begin{Remark} The classification of the (infinitely many) linear types of (3,9)-types was recently carried out in \cite{borovoiClassificationRealTrivectors2022}. Similar methods yield partial information on the (4,8)-case, cf. \cite{leOrbitsRealZgraded2011,leClassificationKformsRn2020}.
\end{Remark}

\section{Flatness of multisymplectic forms}
\label{sec:flat}

In this section we are going to give an overview over various flatness theorems for multisymplectic manifolds in the literature.

\subsection{Helpful tools}

\subsubsection*{Moser trick}

Many of the theorems in this section can be proven using the Moser trick, which was originally applied to volume forms \cite{moserVolumeElementsManifold1965}  and then for proving the symplectic Darboux theorem and its generalization \cite{weinsteinSymplecticManifoldsTheir1971}, cf. Theorem \ref{thm:darbouxclassical} below. We will phrase the trick as the following Lemma:

\begin{Lemma}\label{lem:moser} Let $U\subset \mathbb R^n$, $p\in U$ and $\omega\in \Omega^k(U)$ a multisymplectic form. Let $\alpha$ be a primitive of $\omega-\omega_p$ near $p$, with $\alpha_p=0$. Let us write $\omega_t$ for the time-dependent form $t\omega+(1-t)\omega_p$.\footnote{Here we reinterpret the linear form $\omega_p\in \Lambda^kT_p^*U\cong \Lambda^k\mathbb R^n$ as a constant coefficient form in $\Omega^k(U)\cong C^\infty(U,\Lambda^k\mathbb R^n)$.}
If there exists a time-dependent vector field $X_t$ satisfying the equation:
\begin{align}\label{eq:mosertric}
	\iota_{X_t}\omega_t=\alpha 
\end{align}
near $p$ for $t\in [0,1]$, then $\omega$ is flat near $p$.
\end{Lemma} 
\begin{proof}
By construction, we have $\frac{d\omega_t}{dt}=\omega-\omega_p=d\alpha$. Let $\phi^t$ be the flow of $X_t$. At $p$ we have $(\omega_t)_p=\omega_p$ for all $t$, i.e. $(\omega_t)_p$ is non-degenerate, and thus $\omega_t$ is non-degenerate near $p$. This means that $X_t$ is unique near $p$. Moreover, since $\alpha_p=0$, we have $X_t(p)=0$ for all $t$. This means that the flow of $X_t$ at $p$ exists for all times. In particular, there is a neighborhood of $p$, such that the time-dependent flow of $X_t$ for starting time $t_0=0$ exists until time 1. We write $\phi_t:=\phi_{0,t}$ for this flow. Using (the time-dependent version of) Cartan's formula we can calculate:

\[
\frac{d}{dt}(\phi^t)^*\omega_t=(\phi^t)^*\mathcal L_{X_t}\omega_t +(\phi^t)^*\frac{d\omega_t}{dt}=(\phi^t)^*\left(d\iota_{X_t}\omega_t +\frac{d\omega_t}{dt}\right)=0.
\] 
This means that near $p$ we have $(\phi^1)^*\omega_1=(\phi^0)^*\omega_0=\omega_0$, which is a constant coefficient form. I.e. $\omega=\omega_1$ is flat near $p$.

\end{proof}
\begin{Remark} The Moser trick can make more powerful global statements, and need not be local, nor does it need to intertwine a form with its constant coefficient version. In principle it is applicable to any setting where we have a path of closed 1-forms and appropriate conditions assuring the existence of $X_t$ and its flow. However, here we will only care for the local statement, hence the above formulation.
\end{Remark}

\subsubsection*{$G$-structures}
Another tool, which we are going to need are $G$-structures. Let $M$ be an $n$-dimensional manifold and $G\subset GL(\mathbb R^n)$ a closed Lie subgroup. The bundle of frames $Fr(M)$ is the $ GL(\mathbb R^n)$-principal bundle over $M$, which over $p$ has fiber $Fr(M)_p:=Iso(\mathbb R^n,T_pM)$, where $Iso$ stands for vector space isomorphisms. A $G$-structure on $M$ is then a subspace $P\subset Fr(M)$, which is a $G$-principal bundle with respect to the induced action. We refer to \cite{sternbergLecturesDifferentialGeometry1964} for an introduction to $G$-structures, the statements we will need about them can be summed up by the following:

\begin{Lemma}\label{lem:Gstr} Let $\omega\in\Omega^k(M)$ be of constant linear type $\omega_0\in \Lambda^k ( \mathbb R^n)^*$. Then:
$$
P=\{\phi\in Fr(M)|~\phi^*\omega=\omega_0\}
$$
defines a principal bundle for the group $G=\{g\in  GL(\mathbb R^n)~|~g^*\omega_0=\omega_0\}$. Moreover: 
\begin{itemize}
\item Let $$\omega_0=\sum_{i_1<...<i_k}
\lambda^{i_1,...,i_k}  
 e^{i_1}\wedge ....\wedge e^{i_k}$$
 Then, locally $\omega$ can be written as
 $$\omega=\sum_{i_1<...<i_k}
 \lambda^{i_1,...,i_k}  
 \alpha^{i_1}\wedge ....\wedge \alpha^{i_k}$$
for some  one-forms $\alpha_1,...,\alpha_n$.
\item Let $D_0\subset \mathbb R^n$ be a subspace stabilized by $G$. Then 
$$D=\{\phi(v)|~v\in D_0, \phi\in P\}
$$
defines a smooth distribution on $M$.
\end{itemize}
\end{Lemma} 

\begin{proof}
The fact that $P$ is a principal bundle has been proven in \cite{kobayashiRemarkExistenceGStructure1965} for general tensors of constant linear type. The second statement can be obtained by picking a local section $\sigma:M\supset U\to P$, and setting $(\alpha_i)_p=(\sigma(p)^{-1})^*e^i$. The last statement can be proven by seeing $D$ as the image of the constant rank bundle map $P\times D_0\to TM, (p,v)\mapsto \phi(v)$. 	
\end{proof}

\subsubsection*{Bigrading on the de Rham complex}
In this section we are going to discuss how a decomposition of the (sometimes complexified) tangent bundle induces a bigrading on the de Rham complex. This construction is well-known in the context of almost-complex structures but has been studied in much wider generality in \cite{lychaginNonholonomicFiltrationAlgebraic1994}.

A priori, the de Rham complex $\Omega(M)$ of a manifold has one $\mathbb N_0$-grading and the de Rham differential has degree +1. However, imagine we are given a decomposition of of the tangent into two vector bundles, i.e. $TM=E\oplus F$. Then we have $\Lambda^kT^*M= \bigoplus_{l+r=k}\Lambda^lE^*\otimes \Lambda^rF^*$, and hence:
\begin{align*}
	\Omega^k(M)= \bigoplus_{l+r=k}\Omega^{(l,r)}(M), &&\mathrm{where}&& \Omega^{(l,r)}(M)=\Gamma(\Lambda^lE^*\otimes \Lambda^rF^*).
\end{align*}
With respect to this splitting, the de Rham differential decomposes into multiple components. Since the de Rham differential (and its bidegree) are completely determined on what they do on 1-forms, the only possible bidegrees are $(-1,2),(0,1),(1,0),(2,-1)$. We will denote the corresponding components of $d$ by $L, d^F,d^E, R$. The two 'surprising' components $L$ and $R$ measure exactly the non-involutivities of $E$ resp. $F$.

\begin{Lemma}\label{lem:bigrad} $E$ is involutive if and only if $R=0$.
\end{Lemma}

\begin{proof}Involutivity of $E$ means exactly that for $X,Y\in\Gamma(E)$, we have $[X,Y]\in \Gamma(E)$. Being in $\Gamma(E)$ means exactly lying in the kernel of all $\alpha\in\Omega^{0,1}(M)$. Picking one of these $\alpha$ and $X,Y\in\Gamma(E)$, we have:
$$d\alpha(X,Y)=X(\alpha(Y))-Y(\alpha(X))-\alpha([X,Y])= -\alpha([X,Y]).$$	
Hence $d\alpha$ has no (2,0) component for all $\alpha\in\Omega^{0,1}(M)$ (which is equivalent to $R=0$), if and only if $E$ is involutive.
\end{proof}

\begin{Remark} \label{rem:generalizebigraderham}The above is true much more generally. In the sequel we will make use of two generalizations:
	\begin{itemize}
		\item We can replace $TM$ with any other Lie algebroid, for instance with $TM\otimes \mathbb C$. In this case, an almost-complex structure on $M$ (i.e. $J\in\Gamma(End(TM))$ with $J^2=-id$, induces a splitting of $TM\otimes \mathbb C$, hence a bigrading on the complex-valued de Rham complex. The integrability of the almost-complex structure then is equivalent to the vanishing of $R$ (or equivalently $L$).
		\item We could consider a splitting of $TM$ into more than two components, and would obtain a multigrading instead of a bigrading. 
	\end{itemize}
\end{Remark}

\subsection{Symplectic and volume forms}

Let us start with the most classical cases:

\begin{Theorem}[\cite{moserVolumeElementsManifold1965,weinsteinSymplecticManifoldsTheir1971}]\label{thm:darbouxclassical} Let $\omega$ be a volume form or a symplectic form on an $n$-dimensional manifold, then $\omega$ is flat.
\end{Theorem}

\begin{proof} After restricting to some local coordinates, we only need to show that a time-dependent vector field satisfying Equation \eqref{eq:mosertric} can be found. However, the non-degeneracy of $\omega_t$ at $p$ implies that it is also non-degenerate on an open neighborhood $V$ of $p$ (for all $t\in [0,1]$). This means that the map 
\begin{align*}
\iota_\bullet \omega_t:TV \to \Lambda^{n-1}T^*V && \mathrm{resp.}&& \iota_\bullet \omega_t:TV\to T^*V
\end{align*}
is an injective vector bundle morphism. Since the ranks of $TV$, $T^*V$ and $\Lambda^{n-1}T^*V$ coincide, this means that $\iota_\bullet \omega_t$ is invertible, i.e. that Equation \eqref{eq:mosertric} in Lemma \ref{lem:moser} has a solution for any choice of $\alpha$.
\end{proof}

These two cases are the only ones, where $\iota_\bullet \omega_t$ is a bijection. In the later cases, when we want to use the Moser trick, we will need to pick $\alpha$ with some care, to assure the solvability of equation \eqref{eq:mosertric}.

\subsection{Multicotangent bundles}

The next type of multisymplectic manifolds we look at are multicotangent bundles. These lie at the origin of multisymplectic geometry and their Darboux theorem was investigated in \cite{martinDarbouxTheoremMultisymplectic1988}. We start with an example:

\begin{Example}\label{ex:multangent} Let $Q$ be an $m$-dimensional manifold and $1\leq k\leq m$. Then the following defines a $k$-form $\theta$ on $M=\Lambda^kT^*Q$.
	$$
	\theta_{\eta}(v_1,...,v_k):=\eta(T\pi(v_1),...,T\pi(v_k)),$$
	where $\pi:M\to Q$ is the basepoint projection. If we have coordinates $q_1,...,q_m$ on $Q$, then an $\eta\in M$ can be parametrized as $\eta=\sum_{i_1<...<i_k}p^{i_1,...,i_k}dq_{i_1}\wedge ...\wedge dq_{i_k}$, where $p^{i_1,...,i_k}$ are additional $m\choose k$ coordinates parametrizing the fibers of $\pi$. Then $\omega=d\theta\in \Omega^{k+1}(M)$ is multisymplectic. It's local coordinate expression is:
	$$
	\sum_{i_1<...<i_k}dp^{i_1,...,i_k}\wedge dq_{i_1}\wedge ...\wedge dq_{i_k}.
	$$ 
In the case where $k=1$, this yields the canonical symplectic form on the cotangent bundle of $Q$. When $k=m$, there is only one summand and we obtain a volume form. 
\end{Example}

The forms in the above example have, with respect to a basis $e^i, i\in\{1,...,m\}, f^{J}, J\in \{(i_1,i_2,...,i_k), 1\leq i_1<...<i_k\leq m\}$ of $T^*_pM$, the linear type
\begin{align}\label{eq:multitype}
\omega_0=\sum_{i_1<...<i_k} f^{i_1,...,i_k}\wedge e^{i_1}\wedge ...\wedge e^{i_k}.
\end{align}

Our objective now is to characterize forms which look like this locally, and show a Darboux theorem for them. For this we will use the following invariant:
\begin{align*}
	I_r(\omega):=\{v\in TM ~|~\mathrm{rank}(\iota_v\omega)\leq r\}.
\end{align*}
For a fixed point $p\in M$, $I_r(\omega)_p$ is the set of solutions of a system of (homogenous) real algebraic equations, i.e. an algebraic set. The reason we are interested in $I_r$ is that it allows us to recognize certain directions for forms of linear type \eqref{eq:multitype}:

\begin{Lemma}\label{lem:findingW}
Let $\omega_0$ be as in Equation \eqref{eq:multitype} with $1<k<m$. Then:
\begin{enumerate}
	\item $I_m(\omega_0)=\bigcap_{i=1}^m\ker(e^i)$.
	\item It satisfies \begin{align}\label{eq:conditionmartin}
		\forall v,w\in I_m(\omega)_p: \iota_w\iota_v\omega_p=0.
	\end{align} 
\end{enumerate}
\end{Lemma}

\begin{proof}
Let us start with item 1. For the inclusion $\supset$ we can start with $\xi\in \bigcap_i ker(e^i)$. Then $$\iota_\xi\omega_0=\sum_{i_1<...<i_k} f^{i_1,...,i_k}(\xi)\wedge e^{i_1}\wedge ...\wedge e^{i_k}$$
has rank $\leq m$. Conversely, assume $\xi$ does not lie in the kernel. Without loss of generality $e^1(\xi)=1$. We just have to find $m+1$ vectors $v$ such that $\iota_v\iota_\xi\omega_0$ are linear independent. Let us consider multi-indices $I=1|J=(1,j_1,...,j_{k-1})$, where $J=(j_1,...,j_{k-1})$ is a strictly increasing multi-index of length $k-1$ not containing $1$. There are $m-1\choose k-1$ such indices. We can calculate $\iota_{f_{1|J}}\iota_\xi\omega_0$, where the $f_I, e_i$ is the dual basis of $(f^I,e^i)$. We obtain:	
	\begin{align*}
		\iota_{f_{1|J}}\iota_\xi\omega_0=-\iota_\xi\iota_{f_{1|J}}\omega_0
		=-\iota_\xi\left( e^{1}\wedge e^{j_1}\wedge ... \wedge e^{j_{k-1}}\right)=-e^{j_1}\wedge ... \wedge e^{j_{k-1}} - e^1\wedge \iota_\xi(e^{j_1}\wedge ... \wedge e^{j_{k-1}})\\
		=-e^{j_1}\wedge ... \wedge e^{j_{k-1}} + e^1\wedge (...)
	\end{align*}
	These are all linearly independent among each other (as one can see more easily by considering the expression up to multiples of $e^1$). The worst thing that can happen to us is that there are only ${m-1\choose 2-1}={m-1\choose m-1-1}=m-1$ of those. However, since $m\geq 3$ we can also look at the contractions of $e_2$ (resp. $e_3$) into $\iota_\xi\omega_0$:
	\begin{align*}
		\iota_{e_2}\iota_\xi\omega_0=&\sum_{3\leq i_1<...< i_{k-2}}f^{12|\tilde J}\wedge e^{i_1}\wedge ... \wedge e^{i_{k-2}} + e^1 \wedge (...)\\
		\iota_{e_3}\iota_\xi\omega_0=&\sum_{4\leq i_1<...< i_{k-2}}f^{13|\tilde J}\wedge e^{i_1}\wedge ... \wedge e^{i_{k-2}} + e^1 \wedge (...)+ e^2\wedge (...)
	\end{align*}
	where $(...)$ is an arbitrary expression. These two are linearly independent of all the preceding ones and each other. In total, we have found that the set $\{\iota_{f_{1|J}}\iota_\xi\omega~|~|J|=k-1, 1\not \in J\}\cup \{\iota_{e_2}\iota_\xi\omega,\iota_{e_3}\iota_\xi\omega\}$ is linearly independent, i.e. that the rank of $\iota_\xi\omega$ is strictly larger than $m$.			
	
So, the space $I_m(\omega)_p$ recovers the space spanned by $\{f_I\}$. This space does satisfy condition \eqref{eq:conditionmartin}. 
\end{proof}

\begin{Theorem}[\cite{martinDarbouxTheoremMultisymplectic1988}, cf. also \cite{sevestreLagrangianSubmanifoldsStandard2019}] \label{thm:martin} Let $\omega$ be a multisymplectic $k+1$-form on an ${m\choose k}+m$-dimensional manifold, with $1<k<m$. If at every $p$, the algebraic set $I_m(\omega)_p$ is an $m\choose k$-dimensional vector space such that Equation \eqref{eq:conditionmartin} is satisfied, then:
	\begin{itemize}
		\item $\omega$ is of linear type \eqref{eq:multitype}
		\item The space $I_m(\omega)$ forms a distribution.
	\end{itemize} 
If, moreover, the distribution $I_m(\omega)$ is involutive, then $\omega$ is flat. 
\end{Theorem}
\begin{proof}
	The essence of the proof can be found \cite{martinDarbouxTheoremMultisymplectic1988}, we refer to \cite{sevestreLagrangianSubmanifoldsStandard2019} for a detailed self-contained proof. We note that the formulation in these references do not 'find' $I_m(\omega)$ (which we do by Lemma \ref{lem:findingW}), but rather postulate its existence. We note that for the individual vector spaces   $I_m(\omega)_p$ to assemble to a smooth distribution (which is required for the proof to work), one needs to employ Lemma \ref{lem:Gstr}.
\end{proof}

We refer to \cite{cantrijnGeometryMultisymplecticManifolds1999} for a more global viewpoint on forms of this linear type, in particular containing criteria for a multisymplectic manifold to be globally a multi-cotangent bundle. Moreover, multiple studies have been made of Darboux type theorems for multisymplectic forms which arise in the field-theoretic case, i.e. when replacing $Q$ in Example \ref{ex:multangent} with a fiber bundle and imposing certain horizontality/ verticality conditions on the allowed exterior tensors. (\cite{leonTulczyjewsTriplesLagrangian2003,forgerLagrangianDistributionsConnections2013,forgerMultisymplecticPolysymplecticStructures2013}).

One particularly interesting case occurs when $k=m-1$. Then $\omega$ is an $m$-form in a $2m$-dimensional space. More specifically, it is a binary form in the sense of \cite{turielNformsDimension2n2001} (cf. also \cite{vanzuraSpecialNforms2ndimensional2008}). For these, it was observed in \cite{turielNformsDimension2n2001}, that the involutivity of $I_m\omega$ is automatic when $m\geq 4$. In fact this is true more generally:

\begin{Theorem}\label{thm:martininvolutive}In the setting of Theorem \ref{thm:martin}, if $k>2$ or if $k=2$ and $m\geq 6$, then the involutivity of $I_m(\omega)$ is automatic.	
\end{Theorem}
\begin{proof}
By using Lemma \ref{lem:Gstr}, we can locally find 1-forms $\alpha^I\in\Omega^1(U), I\in \{(i_1,...,i_k)|1\leq i_1<...<i_k\leq m\}$ $\beta^i\in\Omega^1(U)$, $i\in\{1,...,m\}$, such that 
	\begin{align*}
\omega=\sum_{i_1<...<i_k} \alpha^{i_1,...,i_k}\wedge \beta^{i_1}\wedge ...\wedge \beta^{i_k}
\end{align*}
We want to show that $I_m(\omega)=\bigcap_i\ker(\beta_i)$ is involutive. We introduce a bigrading on $\Omega(U)$ with the first grading coming from $\alpha$-directions and the second from $\beta$-directions. Let $R$ be the degree (2,-1) part of the de Rham differential $d$. Following Lemma \ref{lem:bigrad}, for $I_m(\omega)$ to be involutive, we need to show that $R$ is 0. For degree reasons $R(\alpha^I)=0$. Let us write:
\begin{align*}
R(\beta_i)=\sum_{I<J}\lambda^{I,J}_i\alpha^I\wedge \alpha^J
\end{align*}
where we use lexicographic ordering on the multiindeces $I,J$, and $\lambda^{I,J}_i$ are functions on $U$. the closeness of $\omega$ implies that $R(\omega)=0$, i.e. we obtain:
	\begin{align}\nonumber
0&=\sum_{i_1<...<i_k} \alpha^{i_1,...,i_k}\wedge 
\left(\sum_{r=1}^k (-1)^{r} R(\beta^{i_r})\wedge \beta^{i_1}\wedge ...\wedge  
\widehat{\beta^{i_r}}
\wedge... \wedge \beta^{i_k}\right)\\
&=\sum_{i_1<...<i_k}\sum_{r=1}^k\sum_{I<J}\label{eq:complicatedaaarh}
(-1)^{r} \lambda^{I,J}_{i_r}
\left(
 \alpha^{i_1,...,i_k}\wedge 
\alpha^I\wedge \alpha^J\wedge \beta^{i_1}\wedge ...\wedge  
\widehat{\beta^{i_r}}
\wedge... \wedge \beta^{i_k}
\right)
\end{align}
We want to show that this already implies that $\lambda^{A,B}_{i}$ are all identically zero. Let us consider $A<B<C$ multiindeces of length $k$ and $J=\{j_1<...<j_{k-1}\}$ a multiindex of length $k-1$. We want to see, what are the possible terms contributing to $\alpha^A\wedge\alpha^B\wedge \alpha^C\wedge \beta^{j_1}\wedge ...\wedge\beta^{j_{k-1}}$. In principle three types of coefficients can appear:$\pm\lambda^{BC}_{i}$ when $A=J\cup\{i\}$ , $\pm\lambda^{AC}_{j}$ when $B=J\cup\{j\}$ and  $\pm\lambda^{AB}_{l}$ when $C=J\cup\{l\}$.
This means that given $A\neq B, l$, if we can find a multiindex $J$ of length $k-1$ such that $l\not\in J$ and $J\not\subset A$ and $J\not\subset B$, only $\pm\lambda^{AB}_l$ will appear as coefficient in front of $\alpha^A\wedge\alpha^B\wedge \alpha^{J\cup l}\wedge \beta^{j_1}\wedge ...\wedge\beta^{j_{k-1}}$, hence Equation \eqref{eq:complicatedaaarh} will imply $\lambda^{AB}_l=0$. Let us see, when we can find such a $J$:
\begin{itemize}
	\item Assume $A\cap B\not\subset\{l\}$ is non-empty, i.e. let it contain some element $p$. Then pick $J$ to be any index of length $k-1$ in $\{1,...,m\}\backslash\{p,l\}$. It satisfies the desired property.
	\item  Assume $A\cap B\subset \{l\}$. 
	\begin{itemize}
		\item If $k>2$, then just pick $a\in A\backslash\{l\}$,  $b\in B\backslash\{l\}$, and pick $J$ to be a multiindex containing $a$ and $b$, but not $l$. It satisfies the desired property.
		\item If $k=2$, then $A\cup B\cup \{l\}$ has at most 5 elements. Hence for $m\geq 6$, there is some $p\not\in A\cup B\cup \{l\}$. Then we can simply take $J$ to be some multiindex containing $p$ but not $l$. It satisfies the desired property.
	\end{itemize}
\end{itemize}
In total this means that $\lambda^{A,B}_l$ vanish for all $A,B,l$, i.e. that $R=0$, i.e. that $I_m(\omega)$ is involutive.
\end{proof}

When $k=m-1$ the above result recovers a result of \cite{turielNformsDimension2n2001}. The following example shows that involutivity is not automatic when $(m,k)=(3,2)$. The remaining cases $(m,k)=(4,2)$ and $(m,k)=(5,2)$ are to our knowledge unknown.

\begin{Example}\label{ex:nonflat} In \cite{turielNformsDimension2n2001} the following form is presented to show that involutivity is not automatic for $(m,k)=(3,2)$:
\begin{align*}
dy_1 \wedge  dx_2 \wedge  dx_3 + dy_2 \wedge  (dx_1 + y_2dy_3) \wedge  dx_3 + dy_3 \wedge  (dx_1 + y_2dy_3) \wedge  dx_2
\end{align*}	
A similar example has been constructed in \cite{vanzuraOneKindMultisymplectic2001}.
\end{Example}

\subsection{Binary forms}

In this section we are going to discuss Darboux theorems for binary forms, following \cite{turielNformsDimension2n2001} (cf. also \cite{vanzuraSpecialNforms2ndimensional2008}). These forms can be characterized as follows:

\begin{Definition} Let $V$ be a $2m$-dimensional vector space for $m\geq 3$ and $\omega\in\Lambda^mV^*$ be non-degenerate. The form $\omega$ is called binary, if there exists $\omega'\in \Lambda^mV^*\backslash \mathbb R\omega$, such that 
	$$\{\iota_v\omega'~|~v\in V\}\subset  \{\iota_v\omega~|~v\in V\}.$$ 
\end{Definition}
In other words, if we write $E(\tilde \omega):=\mathrm{Image}(v\mapsto \iota_v\tilde \omega)$, a non-degenerate form $\omega$ is called binary  if the vector space
$$
Q=\{\tilde \omega|~E(\tilde\omega)\subset E(\omega)\}
$$
is of dimension $>1$. To any element $\tilde \omega$ of $Q$, we can associate a unique endomorphism $J_{\tilde \omega}$ such that
$
\iota_{J_{\tilde \omega}v}\omega=\iota_v\tilde\omega.
$
The non-degeneracy of $\omega$ implies the uniqueness of $J_{\tilde\omega}$. The existence of the 'binary tensor' $J_{\tilde\omega}$ is probably the origin of calling these forms binary. The map 

$$
Q\to End(V),~\tilde\omega\mapsto J_{\tilde\omega}
$$
is linear and maps $\omega$ to the identity endomorphism. We note that we can move $J_{\tilde\omega}$ around freely inside $\omega$:
$$
\omega(J_{\tilde\omega}v_1,v_2,...,v_m)=\omega(v_1,J_{\tilde\omega}v_2,...,v_m)=...= \omega(v_1,v_2,...,J_{\tilde\omega}v_m).
$$
Moreover, through iteration of the above, one can move around any polynomial expression $P(J_{\tilde\omega})$, inside $\omega$ i.e. 
$$
\omega(P(J_{\tilde\omega})v_1,v_2,...,v_m)=\omega(v_1,P(J_{\tilde\omega})v_2,...,v_m)=...
$$
for any polynomial $P$.

Our goal is to show that a binary form $\omega$ can only be of one of three distinct linear types. We start with the following Lemma:

\begin{Lemma}[\cite{turielNformsDimension2n2001}] \label{lem:prodlin} Let $\omega$ be a binary form and $\omega'\in Q\backslash \mathbb R\omega$ . If $J$ has two or more distinct real eigenvalues, then it is of linear type
\begin{align}\label{eq:prodlin}
e^1\wedge ...\wedge e^m + e^{m+1}\wedge ...\wedge e^{2m}.
\end{align}
We call this binary type \emph{product type}.
\end{Lemma}

\begin{proof}
If $J$ has more than one eigenvalue, then its minimal polynomial $\phi$ decomposes as $\phi=\phi_1\cdot\phi_2$. Then $V=V_1\oplus V_2$, where $V_i=\mathrm{Image}(\phi_i(J))$. We can observe that

$$
\omega(V_1,V_2, ... )=\omega(\phi_1(J)\cdot ,\phi_2(J)\cdot , ... )=\omega(\phi_1\phi_2(J),\cdot, ... )=0.
$$  
This means that $\omega=\omega_1+\omega_2$ with $\omega_1$ zero on $V_2$ and vice-versa. Since $\omega$ is an $m$-form, both of the spaces have to be at least $m$-dimensional, and since $V$ is of dimension $2m$, they are exactly $m$-dimensional. The statement now follows from picking appropriate bases of $V_1$ and $V_2$ and putting them together to a basis of $V$.
\end{proof}

\begin{Remark}\label{rem:prodlin}
We note that for the setting of the previous Lemma, i.e. for $\omega$ of linear type \eqref{eq:prodlin} and $m\geq 3$, any $J$ must be of the shape $a\cdot id_{V_1}\oplus b\cdot id_{V_2}$ for distinct two real eigenvalues $a$ and $b$. The reason for this is that the subspaces $V_1,V_2$ are uniquely characterized by $V_1\cup V_2=I_{m-1}(\omega)=\{v~|~\mathrm{rank}(\iota_v\omega)\leq m-1\}$. Note, that this would fail if we allowed $m=2$.
\end{Remark}

\begin{Lemma}[\cite{turielNformsDimension2n2001}] Let $\omega$ be a binary form and $\omega'\in Q\backslash \mathbb R\omega$ . If $J=J_{\omega'}$ has no real eigenvalues, then it is of linear type
	\begin{align}\label{eq:clxlin}
	\Re((e^1+i f^1)\wedge ...\wedge (e^m+i f^m) ),
	\end{align}
for some basis $e^1,...,e^m,f^1,...,f^m$ of $V^*$ , i.e. it looks like the real part of a complex volume form. We call this binary type \emph{complex type}.
\end{Lemma}
\begin{proof}
If we consider the complexified vector space $V^{\mathbb C}=V\otimes \mathbb C$, then the complex version of Lemma \ref{lem:prodlin} is applicable, and we obtain a decomposition of $\omega^{\mathbb C}$ the complexified version of $\omega$
$$
\omega^{\mathbb C}=\omega_1+\omega_2
$$
with now complex forms $\omega_i$. Since $\omega^{\mathbb C}$ is real, it is stable under complex conjugation, so either $\omega_i=\bar\omega_i$, or $\omega_2=\bar\omega_1$. The former case can not occur, because it would bring us back to the setting of Lemma \ref{lem:prodlin}, which contradicts the non-existence of real eigenvalues (cf. also Remark \ref{rem:prodlin}). In the second situation we obtain the statement of the Lemma upon picking an appropriate complex basis of $V_1$ and using 
$$
\omega^{\mathbb C}=\omega_1+\bar\omega_1=2\Re(\omega_1).
$$
\end{proof}

\begin{Remark}
We will now provide an argument from \cite{vanzuraSpecialNforms2ndimensional2008}, to show that for forms of linear type \eqref{eq:clxlin}, the space $Q$ is also two-dimensional. More precisely the space $J_Q$ of possible $J_{\tilde \omega}$ for $\tilde\omega\in Q$ is given by $\mathbb R\cdot id + \mathbb R\cdot J$, where $J$ is the morphism defined by $e_i\mapsto f_i$, $f_i\mapsto -e_i$, where $(e_i,f_i)$ is the dual basis of $(e^i, f^i)$. 
The ``$\supset$''-inclusion follows from the fact that $J$ behaves like multiplication with $i$. For the other inclusion we first observe that $m\geq 3$ implies that the elements of $J_Q$ commute:
\[
\omega(ABu,v,w, \ldots)=\omega (u,Av,Bw, \ldots)=\omega(BAu,v,w, \ldots).
\]
This implies $AB=BA$ since $\omega$ is non-degenerate. Especially $J_Q\subset End_{\mathbb C}(V)$ (if we use $J$ to define a complex structure on $V$), as every element has to commute with $J$. Moreover,  any element $A\in J_Q$ has to be diagonal as a complex matrix. To see that, we observe that $A(v)$ is always $\mathbb C$-linearly dependent on $v$. We have (again using $m\geq 3$)
\[
\omega(v,Av,x, \ldots)=\omega(v,v, Ax, \ldots)=0,
\]
so $\iota_{v}\iota_{A(v)}\omega=0$ for all $v$.  Now $\omega$  is the real part of a complex volume, so $\iota_{v}\iota_{A(v)}\omega=0$ implies that $v$ is a complex eigenvector of $A$. Moreover, since this is true for any $v$, one can deduce that all eigenvalues are identical.
\end{Remark}

The remaining case is that of a unique real eigenvalue:

\begin{Lemma}[\cite{turielNformsDimension2n2001}] Let $\omega$ be a binary form and $\omega'\in Q\backslash \mathbb R\omega$ . If $J_{\omega'}$ has exactly one real eigenvalue, then it is of linear type
	\begin{align}\label{eq:bin0lin}
		\sum_{i=1}^m e^i \wedge f^1\wedge ... \wedge \widehat{f^i}\wedge ...\wedge f^m
	\end{align}
	for some basis $e^1,...,e^m,f^1,...,f^m$ of $V^*$. We call this binary type \emph{multicotangent type}.
\end{Lemma}
We note that this linear type corresponds to the type described in Equation \eqref{eq:multitype} for $k=m-1$.

\begin{proof} By the same arguments as in the previous Lemmas, $J$ can not have any non-real eigenvalues (because for a real matrix those come in pairs which would give us a totality of $\geq 3$ eigenvalues, which is impossible). Let $\lambda$ be the eigenvalue of $J$. By replacing $J$ with $J-\lambda id$ (and $\omega'$ with $\omega'-\lambda\omega$), we can assume $\lambda=0$. The minimal polynomial of $J$ is $\phi(t)=t^k$ for some $k$. We can now consider instead of $\omega'$ the form $\omega''$, defined by
\[
\omega''(v_1,...,v_m)=\omega(J^{k-1}v_1,...,v_m).
\] 
Hence $J_{\omega''}=J^{k-1}$. This means that $J_{\omega''}^2=0$. Since $\ker(J_{\omega''})=\{v|\iota_v\omega''=0\}$, and a non-zero $m$-form has rank at least $m$, we know that $dim(ker(J_{\omega''}))\leq m$. But since $\mathrm{Image}(J_{\omega''})\subset \ker(J_{\omega''})$, this means that 
$$W=\mathrm{Image}(J_{\omega''})= \ker(J_{\omega''})$$	 
is $m$-dimensional. Pick any complement $U$ of $W$ and a basis $f_1,...,f_m$ of $U$ (and dual basis $\{f^1,...,f^m\}$). Let us consider the map $W\to \Lambda^{m}V^*$, its image lies inside the $m$-dimensional subspace $\Lambda^{m}U^*\subset \Lambda^{m}V^*$, because for two elements $v_1'=J_{\omega''}v_1,v_2'=J_{\omega''}v_1 $ in the image of $J_{\omega''}$, we have 
$$
\iota_{J_{\omega''}v_1}\iota_{J_{\omega''}v_2}\omega=0.
$$
Moreover, since $\omega$ is non-degenerate this map is injective, hence an isomorphism $\varphi:W\to \Lambda^mU^*$. Setting
 $$e_i:=\varphi^{-1}(f^1\wedge ... \wedge \widehat{f^i}\wedge ...\wedge f^m)$$
we obtain the desired basis, such that $\omega$ is of the form \eqref{eq:bin0lin}.
\end{proof}

We can now state the Darboux theorems for binary forms:
\begin{Theorem}[\cite{turielNformsDimension2n2001}]\label{thm:binary} Let $\omega\in \Omega^m(M)$, $m\geq 3$, be a multisymplectic form of constant binary linear type. Then $\omega$ is flat under the following conditions:
	\begin{itemize}
		\item If $\omega$ is of product type \eqref{eq:prodlin}, then the space $I_{m-1}(\omega)$, locally is the union of two distributions. The condition is for one (or equivalently both) of these to be involutive.
		\item If $\omega$ is of complex type \eqref{eq:clxlin}, then it locally induces two almost-complex structures satisfying $\iota_w\iota_{Jv}\omega=\iota_{Jw}\iota_{v}\omega$. The condition is for one (or equivalently both) of these to be integrable.
		\item If $\omega$ is of multicotangent type \eqref{eq:bin0lin}, then it induces a distribution $I_{m-1}(\omega)$  The condition is for this distribution to be involutive.
	\end{itemize}
When $m\geq 4$ these conditions are always satisfied.
\end{Theorem}

\begin{proof}
The necessity of the conditions comes from the fact that constant coefficient forms induce constant coefficient (and hence involutive/integrable) distributions and almost-complex structures. We will prove the three cases separately:
\begin{itemize}
\item In the first case we use Lemma \ref{lem:Gstr} to locally write $\omega$ as
$$
\omega=\omega_\alpha+\omega_\beta=\alpha^1\wedge ...\wedge \alpha^n+ \beta^1\wedge ...\wedge \beta^n
$$
for some coframe $\alpha^i,\beta^i\in \Omega^1(M)$. We have the two distributions $V_1=\bigcap_i \ker(\alpha^i)$ and $V_2=\bigcap_i \ker(\beta^i)$, whose union yields $I_{m-1}(\omega)$. Let us pick the bigrading on differential forms induced by the splitting $TM=V_1\oplus V_2$. With the notations of Lemma \ref{lem:bigrad}, we obtain:
\begin{align*}
0=d\omega=L\omega+d^{V_1}\omega+d^{V_2}\omega+R\omega=&&d^{V_2}\omega_\alpha+&&R\omega_\alpha+ &&L\omega_\beta+&&d^{V_1}\omega_\beta\\
&&(m,1)&&(m-1,2)&&(2,m-1)&&(1,m)
\end{align*}
When $m\geq 4$ all terms in the above equation have different bidegrees, i.e. have to vanish separately. However $R(\omega_\alpha)=0$ implies $R=0$ and $L\omega_\beta=0$ implies $L=0$, showing that both distributions are involutive. When $m=3$, we have one equation with 2 terms $R(\omega_\alpha)+L\omega_\beta=0$. So if one of the distributions is involutive the other will be involutive automatically. 

Once the two distributions are involutive, we can pick coordinates $x_1,...,x_m,y_1,...,y_m$ such that $\partial_{x_i}$ span $V_2$ and $\partial_{y_i}$ span $V_1$. Then we automatically have:
$$
\omega= f(x,y)dx_1\wedge ...\wedge dx_m+ g(x,y)dy_1\wedge ...\wedge dy_m.
$$
The closedness of $\omega$ implies that $f=f(x)$ and $g=g(y)$, and then the flatness follows from applying the Darboux theorem for volume forms (Theorem \ref{thm:darbouxclassical}) for both summands separately.

\item In the second case, we start by using Lemma \ref{lem:Gstr} to locally obtain an almost-complex structure satisfying $\iota_w\iota_{Jv}\omega=\iota_{Jw}\iota_{v}\omega$. The eigenspaces of $J$ yield a splitting of the complexified tangent bundle. Now, the involutivity argument from the previous case carries over upon replacing the de Rham complex with its complexified version $\Omega^{k}_{\mathbb C}(M)=\Gamma(M,\Lambda^k_{\mathbb C}(TM\otimes \mathbb C)^*)$. Using this, we can assure that for $m\geq 4$, the almost-complex structure $J$ is integrable (and in the case $m=3$ we require it). \\
Now we pick complex coordinates on $M$ and obtain $$\omega=f(z,\bar z)dz_1\wedge ...\wedge dz_m + \overline{f(z,\bar z)}d\bar z_1\wedge ...\wedge d\bar z_m .$$
For degree reasons the summands are separately closed, i.e. $f=f(z)$ is holomorphic. Let $F$ be a holomorphic function with $\frac{\partial F}{\partial {z_1}}=f$, then $\tilde z_1=F(z)$, $\tilde z_2=z_2$,..., $\tilde z_m=z_m$ defines holomorphic coordinates, which induce real coordinates with respect to which 
$$\omega = 2 \Re(d\tilde z_1\wedge ... \wedge d\tilde z_m)$$
has constant coefficients.

\item The third case is a special case of Theorems \ref{thm:martin} and \ref{thm:martininvolutive}. 
\end{itemize}	
\end{proof}

We note that in dimension 6, i.e. for $m=3$, all non-degenerate forms are binary, (cf. Appendix \ref{a36}). They have been studied separately in \cite{vanzuraOneKindMultisymplectic2001,panak3formsAlmostComplex2003} and then together in \cite{buresUnifiedTreatmentMultisymplectic2004}. An interesting different description of the complex case can also be found in \cite{hitchinGeometryThreeFormsSix2000}.

\begin{Example}[\cite{turielNformsDimension2n2001}]
We have already seen in Example \ref{ex:nonflat} a binary multisymplectic form not being flat when $m=3$. For the other types they can be given by :
\begin{align*}
&(dx_1+ydy_3)\wedge dx_2\wedge dx_3 + (dx_1-x_2dx_3)\wedge dy_2\wedge dy_3&\mathrm{and}\\
&\Re\left( (dx_1+y_2dx_3+idy_1)\wedge(dx_2+idy_2)\wedge(dx_3+idy_3)\right).
\end{align*}
\end{Example}

\begin{Example}
A more interesting example of a non-flat multisymplectic 3-form of complex type can be constructed on $S^6$ as the restriction of the form\footnote{which is exactly the constant coefficient model of the 8. linear type in Appendix \ref{a37}. We write here $dx^{ijk}$ as a shorthand for $dx_i\wedge dx_j\wedge dx_k$.}
$$
dx^{123} + dx^{145} - dx^{167} + dx^{246} + dx^{257} + dx^{347} - dx^{356}
$$
on $\mathbb R^7$ to the unit sphere. A detailed study in of this form, including a proof of its non-flatness can be found in \cite{wagnerHamiltonianDynamicsGeometry2025}.   
\end{Example}

\begin{Example} For $m=3$ a multisymplectic binary form can change its linear type. A first example of that has been given in \cite{panak3formsAlmostComplex2003}. We will present here a version from \cite{ryvkinMultisymplecticManifoldNot2016,ryvkinInvitationMultisymplecticGeometry2019}:	
$$
dx_1\wedge dx_3\wedge dx_5-dx_1\wedge dx_4\wedge dx_6-dx_2\wedge dx_3\wedge dx_6+x_2dx_2\wedge dx_4\wedge dx_5.
$$	
For $x_2>0$ it has product type, for $x_2=0$ multicotangent type and for $x_2<0$ complex type. To see this, one can utilize the following procedure (cf. \cite{hitchinGeometryThreeFormsSix2000,bryantGeometryAlmostComplex2006}):
\begin{itemize}
	\item Fix a volume form (e.g. $\Omega=dx_1\wedge ...\wedge dx_6)$. 
	\item Compute the unique endomorphism $J$ satisfying $\iota_{Jv}\Omega=(\iota_v\omega)\wedge \omega$. One can check, that it lies in $J_Q\backslash \mathbb Rid$, where $J_Q=\{J_\omega|\omega\in Q\}$ (independently of the binary type $\omega$ has).
	\item Compute the eigenvectors of $J$.
\end{itemize}
\end{Example}
It would be interesting to see, if the above phenomenon can occur for $m\geq 4$.

In case of product type multisymplectic manifolds, a very similar statement as the above theorem holds, even when one allows more factors, i.e. if one looks at a form of the type:
\begin{align}\label{eq:prod}
	e^{1,1}\wedge ...\wedge e^{1,k}+e^{2,1}\wedge ...\wedge e^{2,k}+...+ e^{m,1}\wedge ...\wedge e^{m,k}
\end{align}
on a $k\cdot m$-dimensional manifold with $k\geq 2$ and $m\geq 3$.
\begin{Theorem}[\cite{ryvkinInvitationMultisymplecticGeometry2019}]\label{thm:prod} Let $\omega\in \Omega(M)$ be a multisymplectic form of linear type \eqref{eq:prod}. Then it can be locally written as:
$$\omega=\sum_{i=1}^k\omega_i$$
for decomposable (i.e. rank $m$) forms $\omega_i$. Then $\omega$ is flat if and only if all individual $\omega_i$ are closed.\\

When $m\geq 4$, this is individual closednesses are automatic.
\end{Theorem}

\begin{proof}
Without the '$m\geq 4$'-statement a proof can be found in \cite{ryvkinInvitationMultisymplecticGeometry2019}. It can also be proven by the same argumentation as in Theorem \ref{thm:binary}, cf. also Remark \ref{rem:generalizebigraderham}.
\end{proof}

\subsection{Density-valued symplectic forms}

Next we will treat forms which (in each tangent space) look like the external product of a symplectic form with a volume form. These have recently been examined in \cite{leskiDensityvaluedSymplecticForms2025}. The case where the volume is 1-dimensional has already been studied in \cite{turielClassificationLocale3formes1984} and \cite{vanzuraCharacterizationOneType2004}. The Darboux theorem in this context is the following:

\begin{Theorem}[\cite{leskiDensityvaluedSymplecticForms2025}] \label{thm:dansymp}Let $M$ be $2m+r$-dimensional with $m>1$ and $r\geq 1$. and $\omega$ a multisymplectic $2+r$-form. Assume that the space 
	$$
	F(\omega)_p=\{\alpha\in T_p^*M ~|~\alpha\wedge \omega_p=0\}
	$$
has constant dimension $r$ for all $p\in M$. Then $\omega$ has the following linear type:
\begin{align}\label{eq:lindec}
	e^1\wedge ...\wedge e^r\wedge (f^1\wedge f^2+...+f^{2m-1}\wedge f^{2m}),
\end{align}
where $f^i, i\in\{1,...,2m\}$, $e^j, j\in\{1,...,r\}$ are some basis of $T^*_pM$.\\

If the annihilator of $F(\omega)$ is an involutive distribution, then $\omega$ is flat. This involutivity is automatic, when $m>2$. 
\end{Theorem}

\begin{proof}
	We refer to \cite{leskiDensityvaluedSymplecticForms2025} for a complete proof, but give a short version here.
	Let us start with the linear statement. Let $\omega_0$ a non-degenerate $2+r$-form in a vector space of dimension $(2m+r)$ be given such that $F(\omega_0)$ is $r$-dimensional. Let $e^1,...,e^r$ be a basis of it and $f^1,...,f^{2m}$ be a basis of its complement.
	By definition of $F(\omega_0)$ we know that $$\omega_0=e^1\wedge ... \wedge e^r\wedge \sum_{ij}\lambda_{i,j}f^i\wedge f^j.$$ The non-degeneracy of $\omega_0$ implies that $\sum_{ij}\lambda_{i,j}f^i\wedge f^j$ is non-degenerate (on $\bigcap_i \ker(e^i)$). Then the symplectic basis theorem implies that there is a linear coordinate change such that $\omega_0$ takes the form \eqref{eq:lindec}. Hence, locally on a neighborhood $U$, we can pick a coframe $\alpha^1,...,\alpha^r,\beta^1,...,\beta^{2m}$ such that 
	$$\omega=\omega^\alpha\wedge\omega^\beta=\alpha^1\wedge ...\wedge \alpha^r\wedge(\beta^1\wedge \beta^2+...+\beta^{2m-1}\wedge \beta^{2m}). $$
	We pick the decomposition $L+d^\alpha+d^\beta+R$ of $d$ with respect to this coframe. For the involutivity of the annihilator $\bigcap_i \ker(\alpha^i)$ of $F(\omega)$, we need to show that $R=0$. The form $\omega$ is of degree $(r,2)$, hence $d\omega=0$ implies $R\omega=0$. The latter is equivalent to $R(\omega^\alpha)\wedge \omega^\beta=0$. When $m\geq 3$, Lepage's decomposition (cf. e.g. \cite{libermannSymplecticGeometryAnalytical1987}) implies that $\wedge \omega: \Omega^{r,2}(U)\to \Omega^{r,4}(U)$ is injective\footnote{
		A direct proof can be obtained as follows: It is sufficient to show the injectivity of the map $L:\Lambda^kV^*\to\Lambda^{2n-k}V^*$, $L(\alpha)=\omega^k\wedge \alpha$ for a symplectic vector space $(V,\omega)$ and of dimension $2n$, with $k<n$. Since the domain and codomain of $L$ have the same dimension, we can show surjectivity rather than injectivity. Let $\eta\in \Lambda^{2n-k}V^*$ Since $\omega^n$ is a volume form, there is $\xi=\sum_{i=1}^N\xi^i_1\wedge ...\wedge \xi^i_k$ with $\iota_\xi(\omega^n)=\eta$. For a decomposable $\xi$ (i.e. $N=1$) it can be directly verified that $\iota_\xi(\omega^n)$ is a multiple of $\omega^k$, the statement for general $\xi$ follows by linearity.
		}, i.e. we get $R(\omega^\alpha)=0$, which using $\omega^\alpha=\alpha^1\wedge ...\wedge \alpha^r$ implies thar $R(\alpha^i)=0$ for all $i$, i.e. $R=0$.
	
	So, let us now assume that $F(\omega)$ has involutive annihilator, i.e. we may pick coordinates ($(x_1,...,x_r,y_1,...,y_{2m})$) such that $F(\omega)=\mathrm{span}(dx_1,...,dx_r)$. Then 
	$$
	\omega= dx_1\wedge ...\wedge dx_r \wedge \omega^y, ~~~\omega^y=\sum_{i<j}f(x,y)dy_i\wedge dy_j 
	$$
	We now split $d=d^x+d^y$ and note that $d^y$ alone also satisfies a Poincaré Lemma (cf. e.g. \cite[Chapter 5]{vaismanCohomologyDifferentialForms2016}). Hence $d^y\omega=0$ implies that there is a potential $\eta\in\Omega^{r,1}(U)$ satisfying $d^y\eta=\omega$. Without loss of generality we assume this potential to vanish at $p\in U$. For degree reasons we also have $d\eta=\omega$. Now we can apply Lemma \ref{lem:moser}: The equation $\iota_{X_t}\omega_t=\eta$ will always have a solution (even one where the vector field has only $\partial_{y_j}$ components), hence $\omega$ is flat.
\end{proof}

When $k=1$ this reduces to a theorem proven in \cite{turielClassificationLocale3formes1984} under a slightly more restrictive condition. Therein, the case $(m,r)=(2,1)$ has been studied in great detail, including a full classification of infinitesimally transitive germs of closed forms of this linear type. In particular, the article contains (for $r=1$), the following example showing that the involutivity is not automatic for $m=2$:

\begin{Example} The following $r+2$-form in $\mathbb R^{r+4}$ is not flat:
	$$
	(dx_1\wedge dx_2+dx_3\wedge dx_4)\wedge (dy_1+x_2dx_4)\wedge dy_2\wedge ...\wedge dy_r.
	$$
\end{Example}

\subsection{Multisymplectic forms of degree $n-2$}
Finally, we will take a look at multisymplectic forms of codegree two (i.e. $n-2$ forms in $n$-dimensional manifolds). A Darboux theorem for them has been recently established in \cite{derdzinskiParallelDifferentialForms2025}. It can be formulated as follows:

\begin{Theorem}[\cite{derdzinskiParallelDifferentialForms2025}]\label{thm:codeg} Let $M$ be $n$-dimensional and $\omega$ a multisymplectic form of degree $n-2$. Assume that the space 
	$$
	F(\omega)_p=\{\alpha\in T_p^*M ~|~\alpha\wedge \omega_p=0\}
	$$
	has constant dimension $r$. Then $\omega$ has the following linear type:\footnote{There is one situation, where it is actually one of two types: If $r=0$ and $m\equiv 2\mod 4$, then the $\pm$ below matters and $\omega$ has one of the two linear types, cf. Lemma \ref{n-2-lemma}}
	\begin{align}\label{eq:lindecco}
		(\pm) e^1\wedge ...\wedge e^r\wedge \sum_{i=1}^m f^1\wedge f^2 \wedge ... \wedge \widehat{f^{2i-1}}\wedge \widehat{f^{2i}}\wedge ....\wedge f^{2m-1}\wedge f^{2m} ,
	\end{align}
	where $f^i, i\in\{1,...,2m\}$, $e^j, j\in\{1,...,r\}$ are some basis of $T^*_pM$ and $2m+r=n$. In particular, $\omega$ can be locally written as 
	$$
	\omega=\pm \nu \wedge \eta^{m-1}
	$$
	where $\nu$ is decomposable and satisfies $\{v|~\iota_v\nu=0\}=D$ where $D=F(\omega)^{ann}$ is the annihilator of $F(\omega)$ and $\eta$ is a 2-form non-degenerate when restricted to $D$.
	
	 If the annihilator of $F(\omega)$ is an involutive distribution (which is necessary for flatness), we may pick $\nu$ to be closed. Then $\omega$ is flat, if and only if the $\eta$ corresponding to a closed $\nu$ satisfies $\nu\wedge d\eta=0$.
\end{Theorem}

\begin{proof} The linear algebra part of the statement is a reformulation of Lemma \ref{n-2-lemma}. In a neighborhood $U$, we pick a coframe $\alpha^1,...,\alpha^r,\beta^1,...,\beta^{2m}$ such that
$$\omega=\pm \omega^\alpha\wedge\omega^\beta=\pm \alpha^1\wedge ...\wedge \alpha^r\wedge\left( \sum_{i=1}^m \beta^1\wedge \beta^2 \wedge ... \wedge \widehat{\beta^{2i-1}}\wedge \widehat{\beta^{2i}}\wedge ....\wedge \beta^{2m-1}\wedge \beta^{2m}\right).$$
We observe that $\omega^\beta= \eta^{m-1}$ for $\eta=\frac{1}{\sqrt[m-1]{(m{-}1)!\cdot}}(\beta^1\wedge \beta^2+...+\beta^{2m-1}\wedge \beta^{2m})$. 
When $F(\omega)^{ann}$ is involutive, we can assume that $\alpha^i=dx_i$ for some coordinates $x_i$, by appropriately modifying the $\beta^i$, hence obtain that $\nu= \omega^\alpha = dx_1\wedge ... \wedge dx_r$ is closed. If the additional condition is satisfied, then   $\nu\wedge \eta$ is closed and we can apply Theorem \ref{thm:dansymp} to $\nu\wedge \eta$ and obtain local coordinates such that $\nu\wedge \eta$ takes the form:
$$
dy_1\wedge ... \wedge dy_r \wedge \left( \sum_{i=1}^m dx_{2i-1}\wedge dx_{2i}\right)
$$
With respect to those coordinates also $\omega =\pm \nu\wedge \eta^{m-1}$ has constant coefficients.
\end{proof}

To finish up, let us look at an example, where $\eta$ is not closed (and hence $\omega$ is not flat):

\begin{Example}[cf. Theorem 16.1 in \cite{derdzinskiParallelDifferentialForms2025}] We will consider the case where $r=0$, the other case can be obtained by adding some additional coordinates. Let us consider $\mathbb R^{2m}$, $m\geq 3$ with the form:
$$
\eta= e^{x_1}\left( \sum_{i=1}^{m-1} dx_{2i-1}\wedge dx_{2i} \right) + e^{-(m-2)x_1}dx_{2m-1}\wedge dx_{2m}
$$
We have $d\eta\neq 0$, however $d (\eta^{m-1})=0$.
\end{Example}

\begin{Remark} Non-degenerate two-forms $\eta$ in $2m$-dimensional spaces such that $\eta^{m-1}$ is closed are characterized using the language of effective forms in \cite[Corollary 16.6]{libermannSymplecticGeometryAnalytical1987}.
\end{Remark}

\subsection{A few non-flat examples}

In this subsection we are going to list a few examples of non-flat multisymplectic forms beyond the cases presented in the previous subsections.\\

Let $G$ be a real simple Lie group with Lie algebra $\mathfrak g$ and $\langle\cdot,\cdot\rangle$ its Killing form. Then $\langle [\cdot,\cdot],\cdot\rangle\in\Lambda^3\mathfrak g^*$ is invariant. It hence induces a (biinvariant) closed form $\omega\in \Omega^3(G)$. It is multisymplectic and has constant linear type, however is flat only when $\dim (G)=3$ (cf. \cite[Theorem 4.20]{ryvkinInvitationMultisymplecticGeometry2019}). More general manifolds carrying multisymplectic forms with these linear types have been studied in \cite{leGeometricStructuresAssociated2013}.\\

The stabilizer of the eighth linear $(3,7)$-type $\phi_0$ in Appendix \ref{a37} is the Lie group $G_2$. In particular any $G_2$ structure on a 7-manifold $M$ induces a 3-form $\phi$ of this linear type. On the other side,  the Hodge dual $\psi_0$ of $\phi_0$ with respect to the standard metric of $\mathbb R^7$ is still $G_2$-invariant, i.e. we also obtain a 4-form $\psi$ associated to the $G_2$-structure. We call a $G_2$-structure closed if $d\phi=0$ and coclosed if $d\psi=0$. There exist $G_2$-structures which are closed and not coclosed, these have multisymplectic but non-flat $\phi$, and there are $G_2$-structures which are coclosed but not closed, these have multisymplectic, non-flat $\psi$. We refer to \cite{joyceCompactManifoldsSpecial2000} for a textbook treating $G_2$-structures and \cite{finoRecentResultsClosed2022} for a more recent account. Looking at $G_2$-structures multisymplectically goes back at least to \cite{ibortMultisymplecticGeometryGeneric2001}, a more detailed treatment can be found e.g. in \cite{choRemarksHamiltonianStructures2013}.\\

We expect further such examples to exist for other special holonomy groups: If we take the fifth linear type in Appendix \ref{a37} instead of the eighth, then we obtain a $\tilde G_2$ structure instead of a $G_2$-structure. These have been studied in \cite{leExistenceClosed3forms2008} (cf. also the appendix in \cite{leClassificationKformsRn2020}). Spin(7)-manifolds have been studied as multisymplectic manifolds (with a 4-form in dimension 8) in \cite{kennonSpin7manifoldsMultisymplecticGeometry2021}.

\appendix
\section{Non-degenerate forms in low dimensions}\label{app:threeclass}

In this appendix we list the non-degenerate types of three-forms in dimensions six, seven and eight. We will also note which of them are stable. Moreover, for the dimension 7 case we will include the information on whether the stabilizers of the form contain an element of negative determinant, which we need in Lemma \ref{74}. For a vector space of dimension $n$, let $\{e^1,...,e^n\} $ be a basis of $V^*$. We will use the shorthand notation $e^{ijk}$ for $e^i\wedge e^j\wedge e^k$. 

For the classification in the $(3,6)$-case we refer to \cite{gourewitchLalgebreTrivecteur1935} (cf. also  \cite{capdevielleClassificationFormesTrilineaires1973} and \cite{bryantGeometryAlmostComplex2006}) for the  $(3,7)$-case to \cite{westwickRealTrivectorsRank1981} and for the $(3,8)$-case to \cite{djokovicClassificationTrivectorsEightdimensional1983}. The stability of the forms follows from the stabilizer considerations in the individual articles and has been systematically studied in \cite{leManifoldsAdmittingStable2007}.

\subsection{Non-degenerate (3,6)-types} 
\label{a36}
\begin{multicols}{2}
\begin{enumerate}
	\item $e^{123}+e^{456}$, stable
	\item $e^{123}-e^{156}+e^{246}-e^{345}$, stable
	\item $e^{145}+e^{246}+e^{356}$
\end{enumerate}
\end{multicols}

\subsection{Non-degenerate (3,7)-types}\label{three-in-seven}\label{a37}
In the below table, '(+)' means that the stabilizer of this form only contains elements of positive determinant and '(+-)' means that there is at least one matrix with negative determinant in the stabilizer.

\begin{multicols}{2}
\begin{enumerate}
	\item  $e^{127}+e^{134}+e^{256}$, (+-)
	\item  $e^{125}+e^{127}+e^{147}-e^{237}+e^{346}+e^{347}$, (+)
	\item  $e^{123}+e^{145}+e^{167}$, (+-)
	\item $e^{127}-e^{136}+e^{145}+e^{246}$, (+)\\
	\item $e^{123}-e^{145}+e^{167}+e^{246}+e^{257}+e^{347}-e^{356}$, stable,  (+)
	\item $e^{127}-e^{136}+e^{145}+e^{235}+e^{246}$,   (+)
	\item  $e^{125}+e^{136}+e^{147}+e^{237}-e^{246}+e^{345}$,  (+)
	\item $e^{123}+e^{145}-e^{167}+e^{246}+e^{257}+e^{347}-e^{356}$, stable,  (+)
\end{enumerate}
\end{multicols}

Here are the reasons for the respective markings with (+) and (+-).
\begin{enumerate}
	\item 
	The stabilizer of this form contains an element of negative determinant, given by
	\[
	\begin{pmatrix}
		-1 &  0&  0&  0&  0&  0& 0\\ 
		0 &  1&  0&  0&  0&  0& 0\\ 
		0 &  0&  0&  1&  0&  0& 0\\ 
		0 &  0&  1&  0&  0&  0& 0\\ 
		0 &  0&  0&  0&  1&  0& 0\\ 
		0 &  0&  0&  0&  0&  1& 0\\ 
		0 &  0&  0&  0&  0&  0& -1\\ 
	\end{pmatrix}.
	\]

	\item  	The four connected components of the stabilizer of this form are represented by the identity, the following two matrices and their product (see \cite[Theorem 2.2]{salacMultisymplectic3forms7dimensional2018})
	
	\[
	\begin{pmatrix}
		-1 &  0&  0&  0&  0&  0& 0\\ 
		0 &  1&  0&  0&  0&  0& 0\\ 
		0 &  0&  -1&  0&  0&  0& 0\\ 
		0 &  0&  0&  0&  -1&  0& 0\\ 
		0 &  0&  0&  -1&  0&  0& 0\\ 
		0 &  0&  0&  0&  0&  0& 1\\ 
		0 &  0&  0&  0&  0&  1& 0\\ 
	\end{pmatrix}
	,
	\begin{pmatrix}
		-1 &  0&  0&  0&  0&  0& 0\\ 
		0 &  1&  0&  0&  0&  0& 0\\ 
		0 &  0&  -1&  0&  0&  0& 0\\ 
		0 &  0&  0&  0&  0&  0& 1\\ 
		0 &  0&  0&  0&  0&  1& 0\\ 
		0 &  0&  0&  0&  1&  0& 0\\ 
		0 &  0&  0&  1&  0&  0& 0\\ 
	\end{pmatrix},
	\]
	which both have positive determinant. Hence, the stabilizer only contains elements of positive determinant.
	
	\item 
	It is basically a symplectic form wedged with one extra direction. The stabilizer of this form contains an element of negative determinant, given by
	\[
	\begin{pmatrix}
		1 &  0&  0&  0&  0&  0& 0\\ 
		0 &  0&  1&  0&  0&  0& 0\\ 
		0 &  1&  0&  0&  0&  0& 0\\ 
		0 &  0&  0&  0&  1&  0& 0\\ 
		0 &  0&  0&  1&  0&  0& 0\\ 
		0 &  0&  0&  0&  0&  0& 1\\ 
		0 &  0&  0&  0&  0&  1& 0\\ 
	\end{pmatrix}.
	\]
	\item
	The two connected components of the stabilizer of this form are generated by the identity and the matrix (see \cite[Theorem 2.4]{salacMultisymplectic3forms7dimensional2018})
	
	\[
	\begin{pmatrix}
		-1 &  0&  0&  0&  0&  0& 0\\ 
		0 &  1&  0&  0&  0&  0& 0\\ 
		0 &  0&  -1&  0&  0&  0& 0\\ 
		0 &  0&  0&  1& 0&  0& 0\\ 
		0 &  0&  0&  0&  -1&  0& 0\\ 
		0 &  0&  0&  0&  0&  1& 0\\ 
		0 &  0&  0&  0&  0&  0& -1\\ 
	\end{pmatrix},\]
	which has positive determinant. Hence, the stabilizer only contains elements of positive determinant.
	\item
	The stabilizer of this form is isomorphic to $\tilde G_2$, which is connected. Hence, the stabilizer only contains elements of positive determinant. 
	\item 
	The two connected components of the stabilizer of this form are generated by the identity and the matrix (see \cite[Theorem 2.5]{salacMultisymplectic3forms7dimensional2018})
	
	\[
	\begin{pmatrix}
		0 &  1&  0&  0&  0&  0& 0\\ 
		1 &  0&  0&  0&  0&  0& 0\\ 
		0 &  0&  -1&  0&  0&  0& 0\\ 
		0 &  0&  0&  1& 0&  0& 0\\ 
		0 &  0&  0&  0&  0&  1& 0\\ 
		0 &  0&  0&  0&  1&  0& 0\\ 
		0 &  0&  0&  0&  0&  0& -1\\ 
	\end{pmatrix},\]
	which has positive determinant. Hence, the stabilizer only contains elements of positive determinant.
	\item  
	The stabilizer of this form is connected by \cite[Proposition 10]{buresMultisymplecticFormsDegree2003}. Hence, the stabilizer only contains elements of positive determinant.
	\item
	The stabilizer of this form is isomorphic to $ G_2$, which is connected. Hence, the stabilizer only contains elements of positive determinant.
\end{enumerate}
	\subsection{Non-degenerate (3,8)-types}\label{three-in-eight}\label{a38}

\begin{multicols}{2}
	\begin{enumerate}
		\item $e^{127}+e^{138}+e^{146}+e^{235}$
		\item $e^{128}+e^{137}+e^{146}+e^{236}+e^{245}$
		\item $e^{135}+e^{246}+e^{147}+e^{238}$
		\item $-e^{135}+e^{146}+e^{236}+e^{245}+e^{127}+e^{348}$
		\item $e^{138}+e^{147}+e^{156}+e^{235}+e^{246}$
		\item $e^{128}+e^{137}+e^{146}+e^{247}+e^{256}+e^{345}$
		\item $e^{156}+e^{178}+e^{234}$
		\item $e^{158}+e^{167}+e^{234}+e^{256}$
		\item $e^{148}+e^{157}+e^{236}+e^{245}+e^{347}$
		\item $e^{134}+e^{234}+e^{156}+e^{278}$
		\item $e^{135}-e^{245}+e^{146}+e^{236}+e^{678}$
		\item $e^{137}+e^{237}+e^{256}+e^{148}+e^{345}$
		\item $e^{135}+e^{245}+e^{146}-e^{236}+e^{678}+e^{127}$
		\item $e^{138}+e^{147}+e^{245}+e^{267}+e^{356}$
		\item $-e^{135}+e^{146}+e^{236}+e^{245}+e^{137}+e^{247}+e^{568}$
		\item $-e^{135}+e^{146}+e^{236}+e^{245}+e^{127}+e^{347}+e^{568}$
		\item $e^{128}+e^{147}+e^{236}+e^{257}+e^{358}+e^{456}$
		\item $-e^{135}+e^{146}+e^{236}+e^{245}+e^{137}+e^{247}+e^{128}-e^{568}$
		\item $e^{124}+e^{134}+e^{256}+e^{378}+e^{157}+e^{468}$, stable
		\item $e^{135}+e^{245}+e^{146}-e^{236}+e^{127}+e^{348}+e^{678}$, stable
		\item $e^{135}-e^{146}+e^{236}+e^{245}+e^{347}+e^{568}+e^{127}+e^{128}$, stable
	\end{enumerate}

	\end{multicols}

\newpage	

\printbibliography

\end{document}